\documentclass{amsart}

\usepackage{amsthm,amsmath,amssymb,graphicx}
\usepackage[usenames]{color}

\title
[Spacelike CMC $1$ surfaces with elliptic ends in $\mathbb{S}^3_1$]
{Spacelike CMC $1$ surfaces with elliptic ends \\ in de Sitter $3$-Space}
\author{Shoichi FUJIMORI}
\date{July 14, 2005}
\address{Department of Mathematics, Kobe University, Kobe 657-8501, Japan}
\email{fujimori@math.kobe-u.ac.jp}
\dedicatory{Dedicated to Professor Takeshi Sasaki 
            on the occasion of his sixtieth birthday}
\def\transpose#1{\mathord{\mathopen{{\vphantom{#1}}^t}#1}}

\DeclareFontFamily{U}{rsfs}{\skewchar\font"7F}
\DeclareFontShape{U}{rsfs}{m}{n}{
	<-6> rsfs5
	<6-8> rsfs7
	<8-> rsfs10
	}{}
\DeclareMathAlphabet{\mathscr}{U}{rsfs}{m}{n}
 \newtheorem{theorem}{Theorem}[section]
 \newtheorem*{theorem*}{Theorem}
 \newtheorem{proposition}[theorem]{Proposition}
 \newtheorem{corollary}[theorem]{Corollary}
 \newtheorem{lemma}[theorem]{Lemma}

\theoremstyle{definition}
 \newtheorem{definition}[theorem]{Definition}
 \newtheorem{example}[theorem]{Example}
\theoremstyle{remark}
 \newtheorem{remark}[theorem]{Remark}
 \newtheorem*{remark*}{Remark}

\numberwithin{equation}{section}

\begin{document} 

\begin{abstract}
We show that an Osserman-type inequality holds for spacelike surfaces of 
constant mean curvature $1$ with singularities and with elliptic ends 
in de Sitter $3$-space. 
An immersed end of a constant mean curvature $1$ surface is an 
``elliptic end'' if the monodromy representation at the end is diagonalizable 
with eigenvalues in the unit circle.  
We also give a necessary and sufficient condition for equality in the 
inequality to hold, and in the process of doing this we derive a condition for 
determining when elliptic ends are embedded. 
\end{abstract}

\maketitle

%
\section*{Introduction} 

It is known that there is a representation formula, using holomorphic null 
immersions into $SL(2,\mathbb{C})$, for spacelike constant mean curvature 
(CMC) $1$ immersions in de Sitter $3$-space $\mathbb{S}^3_1$ \cite{AA}. 
Although this formula is very similar to a representation formula for CMC $1$ 
immersions in hyperbolic $3$-space $\mathbb{H}^3$ (the so-called Bryant 
representation formula \cite{B, UY1}), and the global properties of CMC $1$ 
immersions in $\mathbb{H}^3$ have been investigated 
\cite{CHR, RUY4, RUY5, UY1, UY2, UY3, Yu}, 
global properties and singularities of spacelike CMC $1$ immersions in 
$\mathbb{S}^3_1$ are not yet well understood. 
One of the biggest reasons for this is that the only complete spacelike 
CMC $1$ immersion in $\mathbb{S}^3_1$ is the flat totally umbilic 
immersion \cite{Ak, R}. 
This situation is somewhat parallel to the relation between minimal immersions 
in Euclicean $3$-space $\mathbb{R}^3$ and spacelike maximal immersions in 
Lorentz $3$-space $\mathbb{R}^3_1$, the $\mathbb{R}^3$ case having been 
extensively studied while the $\mathbb{R}^3_1$ case still being not well 
understood. 
Recently, Umehara and Yamada defined spacelike maximal surfaces with certain 
kinds of ``admissible'' singularities, 
and named them ``maxfaces'' \cite{UYmax}.  They then showed that maxfaces are 
rich objects with respect to global geometry. 
So in Section \ref{sec:CMC1faces}, we introduce admissible singularities of 
spacelike CMC $1$ surfaces in $\mathbb{S}^3_1$ and name these surfaces with 
admissible singularities ``CMC $1$ faces''.

For CMC $1$ immersions in $\mathbb{H}^3$, the monodromy representation at each 
end is always diagonalizable with eigenvalues in 
$\mathbb{S}^1=\{z\in\mathbb{C}\,|\,|z|=1\}$, but this is not true for CMC $1$ 
faces in general. 
However, when the monodromy representation at each end of a CMC $1$ face is 
diagonalizable with eigenvalues in $\mathbb{S}^1$, 
we can directly apply many results for CMC $1$ immersions in $\mathbb{H}^3$ to 
CMC $1$ faces. 
So in Section \ref{sec:unitaryends}, we give a definition of 
``elliptic ends'' of CMC $1$ faces. 
In Section \ref{sec:embeddedness}, 
we give a necessary and sufficient condition for elliptic ends to be embedded. 

The total curvature of complete minimal immersions in $\mathbb{R}^3$ of 
finite total curvature satisfies 
the Osserman inequality \cite[Theorem 3.2]{O}. 
Furthermore, the condition for equality was given in \cite[Theorem 4]{JM}. 
This inequality is a stronger version of the Cohn-Vossen inequality. 
For a minimal immersion, the total curvature is equal to the degree of 
its Gauss map, multiplied by $-4\pi$. 
So the Osserman inequality can be viewed as an inequality about the degree of 
the Gauss map. 
In the case of a CMC $1$ immersion in $\mathbb{H}^3$, the total curvature 
never satisfies equality of the Cohn-Vossen inequality 
\cite[Theorem 4.3]{UY1}, and the Osserman inequality does not hold in general. 
However, using the hyperbolic Gauss map instead of the total curvature, 
an Osserman type inequality holds for CMC $1$ immersions in $\mathbb{H}^3$ 
\cite{UY2}. 
Also, Umehara and Yamada showed that the Osserman inequality holds for 
maxfaces \cite{UYmax}. 
In Section \ref{sec:Ossermanineq}, we give the central result of this paper, 
that the Osserman inequality holds for complete CMC $1$ faces of finite type 
with elliptic ends (Theorem \ref{th:main}). 
{\it 
We remark that the assumptions of finite type and ellipticity of the ends 
can actually be removed, because, in fact, any complete CMC $1$ face must be 
of finite type and each end must be elliptic. 
This deep result will be shown in the forthcoming paper \cite{FRUYY}.
} 

Lee and Yang were the first to construct a numerous collection of 
examples of CMC $1$ faces \cite{LY}. 
Furthermore, they constructed complete irreducible CMC $1$ faces with three 
elliptic ends, by using hypergeometric functions \cite{LY}.  
We will also give numerous examples here in Section \ref{sec:examples}, 
by using a method for transferring CMC $1$ immersions in $\mathbb{H}^3$ to 
CMC $1$ faces in $\mathbb{S}^3_1$.  
Applying this method to examples in \cite{MU, RUY4, RUY5}, we give many 
examples of complete reducible CMC $1$ faces of finite type 
with elliptic ends. 

The author would like to thank Professors Masaaki Umehara, Kotaro Yamada and 
Wayne Rossman for their valuable comments and suggestions, 
as well as Professor Jun-ichi Inoguchi for initially suggesting the 
investigation of an Osserman-type inequality for CMC $1$ faces. 

\section{CMC $1$ faces} 
\label{sec:CMC1faces} 

Let $\mathbb{R}^4_1$ be the $4$-dimensional Lorentz space with the Lorentz 
metric 
\[ 
\langle (x_0,x_1,x_2,x_3),(y_0,y_1,y_2,y_3)\rangle =
-x_0y_0+x_1y_1+x_2y_2+x_3y_3. 
\]
Then de Sitter $3$-space is 
\[
\mathbb{S}^3_1=\mathbb{S}^3_1(1)=
\{(x_0,x_1,x_2,x_3)\in\mathbb{R}^4_1 \, | \, -x_0^2+x_1^2+x_2^2+x_3^2=1\},
\]
with metric induced from $\mathbb{R}^4_1$. 
$\mathbb{S}^3_1$ is a simply-connected 
$3$-dimensional Lorentzian manifold with constant sectional curvature $1$. 
We can consider $\mathbb{R}^4_1$ to be the $2 \times 2$ self-adjoint matrices 
($X^* = X$, where $X^*=\transpose{\overline{X}}$ and $\transpose{X}$ denotes 
the transpose of $X$), by the identification 
\[
\mathbb{R}^4_1\ni X=(x_0,x_1,x_2,x_3)\leftrightarrow 
X=\sum_{k=0}^3 x_k e_k 
 =\begin{pmatrix} x_0+x_3 & x_1+i x_2 \\ 
                  x_1-i x_2 & x_0-x_3 \end{pmatrix},
\]
where
\[
e_0=\begin{pmatrix}1&0\\0&1\end{pmatrix},\quad
e_1=\begin{pmatrix}0&1\\1&0\end{pmatrix},\quad
e_2=\begin{pmatrix}0&i\\-i&0\end{pmatrix},\quad
e_3=\begin{pmatrix}1&0\\0&-1\end{pmatrix}.
\]
Then $\mathbb{S}^3_1$ is 
\[ 
\mathbb{S}^3_1=\{X\,|\,X^*=X\,,\det X=-1\}
              =\{Fe_3F^*\,|\,F\in SL(2,\mathbb{C})\} 
\] 
with the metric
\[
 \langle X,Y\rangle 
=-\frac{1}{2}\mathrm{trace}\left(Xe_2(\transpose{Y})e_2\right) . 
\]
In particular, $\langle X,X\rangle =-\det X$. 
An immersion in $\mathbb{S}^3_1$ is called {\em spacelike} if the induced 
metric on the immersed surface is positive definite.  

Aiyama and Akutagawa \cite{AA} showed the following representation formula for 
simply-connected spacelike CMC $1$ immersions. 

\begin{theorem}[The representation of Aiyama-Akutagawa]\label{th:AA} 
Let $D$ be a simply-connected domain in $\mathbb{C}$ with a base point 
$z_0\in D$. Let 
$g:D\to (\mathbb{C}\cup\{\infty\})\setminus\{z\in\mathbb{C}\,|\,|z|\le 1\}$ be 
a meromorphic function and $\omega$ a holomorphic $1$-form on $D$ such that 
\begin{equation}\label{eq:dshat^2}
d\hat s^2=(1+|g|^2)^2\omega\bar\omega
\end{equation}
is a Riemannian metric on $D$. 
Choose the holomorphic immersion $F=(F_{jk}):D\to SL(2,\mathbb{C})$ 
so that $F(z_0)=e_0$ and $F$ satisfies 
\begin{equation}\label{eq:F^-1dF}
F^{-1}dF=\begin{pmatrix}g&-g^2\\1&-g\end{pmatrix}\omega .
\end{equation}
Then $f:D\to\mathbb{S}^3_1$ defined by 
\begin{equation}\label{eq:f=Fe_3F^*}
f=Fe_3F^*
\end{equation}
is a conformal spacelike CMC $1$ immersion. The induced metric 
$ds^2=f^*(ds^2_{\mathbb{S}^3_1})$ on $D$, the second fundamental form $h$, 
and the hyperbolic Gauss map $G$ of $f$ are given as follows:
\begin{equation}\label{eq:ds^2-h-G}
ds^2=(1-|g|^2)^2\omega\bar\omega ,\quad 
   h=Q+\overline{Q}+ds^2,\quad 
   G=\frac{dF_{11}}{dF_{21}}=\frac{dF_{12}}{dF_{22}},
\end{equation}
where $Q=\omega dg$ is the Hopf differential of $f$. 

Conversely, any simply-connected spacelike CMC $1$ immersion can be 
represented this way. 
\end{theorem}

\begin{remark}\label{rm:AArep}
We make the following remarks about Theorem \ref{th:AA}:
\begin{enumerate}
\item Following the terminology of Umehara and Yamada, $g$ is called the 
{\em secondary} Gauss map. The pair $(g,\omega)$ is called the 
{\em Weierstrass data}. 
\item The oriented unit normal vector field $N$ of $f$ is given as 
\[
N=\frac{1}{|g|^2-1}
  F\begin{pmatrix}|g|^2+1&2g\\2\bar g&|g|^2+1\end{pmatrix}F^*,
\]
which is a future pointing timelike vector field (See \cite{KY}). 
\item\label{rm:hypG} 
The hyperbolic Gauss map has the following geometric meaning: 
Let $\mathbb{S}^2_\infty$ be the future pointing ideal boundary of 
$\mathbb{S}^3_1$. Then $\mathbb{S}^2_\infty$ is identified with the Riemann 
sphere $\mathbb{C}\cup\{\infty\}$ in the standard way. 
Let $\gamma_z$ be the geodesic ray starting at $f(z)$ in $\mathbb{S}^3_1$ 
with the initial velocity vector $N(z)$, 
the oriented unit normal vector of $f(D)$ at $f(z)$. 
Then $G(z)$ is the point $\mathbb{S}^2_\infty$ determined by the asymptotic 
class of $\gamma_z$.  See \cite{B, UY1, UY3}.
\item\label{rm:FF*} 
For the holomorphic immersion $F$ satisfying Equation \eqref{eq:F^-1dF}, 
$\hat f:=FF^*:D\to\mathbb{H}^3$ is a conformal CMC $1$ immersion with first 
fundamental form $\hat f^*(ds^2_{\mathbb{H}^3})=d\hat s^2$ as in 
Equation \eqref{eq:dshat^2}, and with second fundamental form 
$\hat h=-\hat Q-\overline{\hat Q}+d\hat s^2$. 
The Hopf differential $\hat Q$ and the hyperbolic Gauss map $\hat G$ of 
$\hat f$ are the same as $Q$ and $G$, 
since $f$ and $\hat f$ are both constructed from the same $F$, and $F$ 
determines $Q$ and $G$, 
by Equations \eqref{eq:F^-1dF} and \eqref{eq:ds^2-h-G}. 
\item\label{rm:2Q=Sg-SG} 
By Equation (2.6) in \cite{UY1}, $G$ and $g$ and $Q$ have the following 
relation: 
\[
2Q=S(g)-S(G),
\]
where $S(g)=S_z(g)dz^2$ and
\[
S_z(g)=\left(\frac{g''}{g'}\right)'-\frac{1}{2}\left(\frac{g''}{g'}\right)^2
\qquad ({}'=d/dz)
\]
is the Schwarzian derivative of $g$. 
\end{enumerate}
\end{remark}

Our first task is to extend the above theorem to non-simply-connected CMC $1$ 
surfaces with singularities, along the same lines as in \cite{UYmax}.  
To do this, we first define admissible singularities. 

\begin{definition}
Let $M$ be an oriented $2$-manifold.  
A smooth (that is, $C^\infty$) map $f:M\to\mathbb{S}^3_1$ is called a 
{\em CMC $1$ map} if there exists an open dense subset $W\subset M$ such that 
$f|_W$ is a spacelike CMC $1$ immersion. 
A point $p\in M$ is called a {\em singular point} of $f$ if the induced metric 
$ds^2$ is degenerate at $p$. 
\end{definition}

\begin{definition}
Let $f:M\to\mathbb{S}^3_1$ be a CMC $1$ map and $W\subset M$ an open dense 
subset such that $f|_W$ is a CMC $1$ immersion. 
A singular point $p\in M\setminus W$ is called an 
{\em admissible} singular point if 
\begin{enumerate}
\item there exists a $C^1$-differentiable function 
$\beta :U\cap W\to\mathbb{R}^+$, 
where $U$ is a neighborhood of $p$, such that $\beta ds^2$ extends to a 
$C^1$-differentiable Riemannian metric on $U$, and 
\item $df(p)\ne 0$, that is, $df$ has rank $1$ at $p$. 
\end{enumerate}
We call a CMC $1$ map $f$ a {\em CMC $1$ face} if each singular point is 
admissible.
\end{definition}

To extend Theorem \ref{th:AA} to CMC $1$ faces that are not simply-connected, 
we prepare two propositions.  First, we prove the following proposition:

\begin{proposition}\label{pr:UYmax2.3}
Let $M$ be an oriented %
$2$-manifold and $f:M\to\mathbb{S}^3_1$ a CMC $1$ face, where $W\subset M$ an 
open dense subset such that $f|_W$ is a CMC $1$ immersion. 
Then there exists a unique complex structure $J$ on $M$ such that 
\begin{enumerate}
\item\label{pr:UYmax2.3(1)} 
$f|_W$ is conformal with respect to $J$, and
\item\label{pr:UYmax2.3(2)} 
there exists an immersion $F:\widetilde{M}\to SL(2,\mathbb{C})$ which is 
holomorphic with respect to $J$ such that 
\[
\det (dF)=0\quad\text{and}\quad f\circ\varrho =Fe_3F^*, 
\]
where $\varrho :\widetilde{M}\to M$ is the universal cover of $M$. 
\end{enumerate}
This $F$ is called a {\em holomorphic null lift} of $f$. 
\end{proposition}

\begin{remark}\label{rm:SU11ambiguity}
The holomorphic null lift $F$ of $f$ is unique up to right multiplication by 
a constant matrix in $SU(1,1)$. 
See also Remark \ref{rm:FtoFB} below. 
\end{remark}

\begin{proof}[Proof of Proposition \ref{pr:UYmax2.3}]
Since the induced metric $ds^2$ gives a Riemannian metric on $W$, it induces a 
complex structure $J_0$ on $W$. 
Let $p$ be an admissible singular point of $f$ and $U$ a local 
simply-connected neighborhood of $p$.  Then by definition, there exists a 
$C^1$-differentiable function $\beta :U\cap W\to\mathbb{R}^+$ such that 
$\beta ds^2$ extends to a $C^1$-differentiable Riemannian metric on $U$.  
Then there exists a positively oriented orthonormal frame field 
$\{\boldsymbol{v}_1,\boldsymbol{v}_2\}$ with respect to $\beta ds^2$ which is 
$C^1$-differentiable on $U$.  
Using this, we can define a $C^1$-differentiable almost complex structure $J$ 
on $U$ such that 
\begin{equation}\label{eq:almostCS}
J(\boldsymbol{v}_1)= \boldsymbol{v}_2\quad\text{and}\quad
J(\boldsymbol{v}_2)=-\boldsymbol{v}_1.
\end{equation}
Since $ds^2$ is conformal to $\beta ds^2$ on $U\cap W$, 
$J$ is compatible with $J_0$ on $U\cap W$. 
There exists a $C^1$-differentiable decomposition 
\begin{equation}\label{eq:decomp}
\Gamma\left(T^*\!M^\mathbb{C}\!\otimes\mathfrak{sl}(2,\mathbb{C})\right)=
\Gamma\left(T^*\!M^{(1,0)}   \!\otimes\mathfrak{sl}(2,\mathbb{C})\right)\oplus
\Gamma\left(T^*\!M^{(0,1)}   \!\otimes\mathfrak{sl}(2,\mathbb{C})\right)
\end{equation}
with respect to $J$, where $\Gamma (E)$ denotes the sections of a vector 
bundle $E$ on $U$.  Since $f$ is smooth, 
$df\cdot f^{-1}$ is a smooth $\mathfrak{sl}(2,\mathbb{C})$-valued $1$-form. 
We can take the $(1,0)$-part $\zeta$ of $df\cdot f^{-1}$ with respect to this 
decomposition.  
Then $\zeta$ is a $C^1$-differentiable $\mathfrak{sl}(2,\mathbb{C})$-valued 
$1$-form which is holomorphic on $U\cap W$ with respect to the equivalent 
complex structures $J_0$ and $J$ 
(which follows from the fact that $f|_W$ is a CMC $1$ immersion, 
 so the hyperbolic Gauss map $G$ of $f$ is holomorphic on $W$, 
 which is equivalent to the holomorphicity of $\zeta$ with respect to $J$ and 
 $J_0$ on $U\cap W$).  
Hence $d\zeta\equiv 0$ on $U\cap W$.  
Moreover, since $W$ is an open dense subset and $\zeta$ is 
$C^1$-differentiable on $U$, $d\zeta\equiv 0$ on $U$.  
Similarly, $\zeta\wedge\zeta\equiv 0$ on $U$.  In particular, 
\begin{equation}\label{eq:int_cond}
d\zeta +\zeta\wedge\zeta = 0. 
\end{equation}

As $U$ is simply-connected, the existence of a $C^1$-differentiable map 
$F_U=(F_{jk})_{j,k=1,2}:U\to SL(2,\mathbb{C})$ such that 
$dF_U\cdot F_U^{-1}=\zeta$ is equivalent to the condition \eqref{eq:int_cond}. 
Hence such an $F$ exists. 

Note that since $f$ takes Hermitian matrix values, 
we have $df=\zeta f+(\zeta f)^*$. 
So $df(p)\ne 0$ (that is, $p$ is an admissible singularity) implies 
$\zeta\ne 0$. 
Then at least one entry $dF_{jk}$ of $dF_U$ does not vanish at $p$. 
Using this $F_{jk}$, we define the function 
$z=F_{jk}:U\to\mathbb{C}=\mathbb{R}^2$. 
Then, $z$ gives a coordinate system on $U$. 
Since $z=F_{jk}$ is a holomorphic function on $U\cap W$, it gives a complex 
analytic coordinate around $p$ compatible with respect to that of $U\cap W$. 
The other entries of $F_U$ are holomorphic functions with respect to $z$ on 
$U\cap W$ and are $C^1$-differentials on $U$, so each entry of $F_U$ is 
holomorphic with respect to $z$ on $U$, by the Cauchy-Riemann equations. 
Since $p$ is an arbitrary fixed admissible singularity, the complex structure 
of $W$ extends across each singular point $p$. 

This complex structure can be seen to be well-defined at singular points as 
follows: 
Let $p'\in M\setminus W$ be another singular point and $U'$ a neighborhood of 
$p'$ so that $U\cap U'\ne\emptyset$. 
Then by the same argument as above, there exists a $C^1$-differentiable 
almost complex structure $J'$ on $U'$ and $C^1$-differentiabe map 
$F_{U'}'=(F'_{j'k'})_{j',k'=1,2}:U'\to SL(2,\mathbb{C})$ such that 
$dF_{U'}'\cdot {F_{U'}'}^{-1}$ is the $(1,0)$-part of $df\cdot f^{-1}$ with 
respect to Equation \eqref{eq:decomp}.  
Define $z'=F'_{j'k'}$ so that $dF'_{j'k'}\ne 0$. 
Then by uniqueness of ordinary differential equations, $F_U=F_{U'}'A$ for 
some constant matrix $A$.  
So $z$ and $z'$ are linearly related, 
and hence they are holomorphically related. 
Also, because $dz$ and $dz'$ are nonzero, we have $dz/dz'\ne 0$ 
on $U\cap U'$. 

For local coordinates $z$ on $M$ compatible with $J$, 
$\partial f\cdot f^{-1}:=(f_zdz)\cdot f^{-1}$ 
(which is equal to $\zeta$) 
is holomorphic on $M$ and there exists a holomorphic map 
$F:\widetilde{M}\to SL(2,\mathbb{C})$ such that 
\begin{equation}\label{eq:dellff^-1}
dF\cdot F^{-1}=\partial f\cdot f^{-1}.
\end{equation}
Since $\partial f\cdot f^{-1}\ne 0$, also $dF\ne 0$, 
and hence $F$ is an immersion. 
Also, since $f$ is conformal, 
$0=\langle\partial f,\partial f\rangle =-\det (\partial f)$. 
Thus $\det (dF)=0$. 

Finally, we set $\hat f=Fe_3F^*$, defined on $\widetilde{M}$.  
We consider some simply-connected region $V\subset W$. 
By Theorem \ref{th:AA}, there exists a holomorphic null lift $\hat F$ of $f$, 
\begin{equation}\label{eq:hatFe3hatF*}
f=\hat Fe_3\hat F^*, 
\end{equation}
defined on that same $V$. 
Then by Equations \eqref{eq:dellff^-1} and \eqref{eq:hatFe3hatF*}, we have 
\[
d\hat F\cdot\hat F^{-1}=dF\cdot F^{-1},
\]
and hence $\hat F=FB$ for some constant $B\in SL(2,\mathbb{C})$. 
We are free to choose the solution $F$ of Equation \eqref{eq:dellff^-1} 
so that $B=e_0$, that is, $\hat F=F$, so $f=\hat f$ on $V$. 
By the holomorphicity of $F$, $\hat f$ is real analytic on $\widetilde{M}$.  
Also, $f\circ\varrho$ is real analytic on $\widetilde{M}$, 
by Equation \eqref{eq:dellff^-1} and the holomorphicity of $F$. 
Therefore $f\circ\varrho =\hat f$ on $\widetilde{M}$, 
proving the proposition.  
\end{proof}

By Proposition \ref{pr:UYmax2.3}, the $2$-manifold $M$ on which a CMC $1$ face 
$f:M\to\mathbb{S}^3_1$ is defined always has a complex structure.  
So throughout this paper, we will treat $M$ as a Riemann surface with a 
complex structure induced as in Proposition \ref{pr:UYmax2.3}.  

The next proposition is the converse to Proposition \ref{pr:UYmax2.3}: 

\begin{proposition}\label{pr:UYmax2.4}
Let $M$ be a Riemann surface and $F:M\to SL(2,\mathbb{C})$ a 
holomorphic null immersion.  Assume the symmetric $(0,2)$-tensor 
\begin{equation}\label{eq:nondegeneracy}
\det[d(Fe_3F^*)]
\end{equation}
is not identically zero. 
Then $f=Fe_3F^*:M\to\mathbb{S}^3_1$ is a CMC $1$ face, 
and $p\in M$ is a singular point of $f$ if and only if 
$\det[d(Fe_3F^*)]_p=0$.  Moreover, 
\begin{equation}\label{eq:d(FF^*)}
-\det[d(FF^*)]\quad\text{is positive definite}
\end{equation}
on $M$. 
\end{proposition}

\begin{proof}
Since \eqref{eq:nondegeneracy} is not identically zero, the set 
\[
W:=\{p\in M|\det[d(Fe_3F^*)]_p\ne 0\}
\]
is open and dense in $M$.  
Since $F^{-1}dF$ is a $\mathfrak{sl}(2,\mathbb{C})$-valued $1$-form, 
there exist holomorphic $1$-forms $a_1$, $a_2$ and $a_3$ such that 
\[
F^{-1}dF=\begin{pmatrix}a_1&a_2\\a_3&-a_1\end{pmatrix}.
\]
Since $F$ is a null immersion, that is, 
$\text{rank}(dF)=1$, $a_j$ ($j=1,2,3$) satisfy 
\begin{equation}\label{eq:a1a2a3}
a_1^2+a_2a_3=0 \qquad\text{and}\qquad |a_1|^2+|a_2|^2+|a_3|^2>0.
\end{equation}
Since 
\begin{align*}
d(Fe_3F^*)
&= F\left( F^{-1}dFe_3+(F^{-1}dFe_3)^*\right) F^* \\
&= F\begin{pmatrix}a_1+\overline{a_1}&-a_2+\overline{a_3}\\
                   a_3-\overline{a_2}&a_1+\overline{a_1}\end{pmatrix}F^*,
\end{align*}
we have 
\begin{align*}
-\det[d(Fe_3F^*)]
&= -2|a_1|^2+|a_2|^2+|a_3|^2 \\ 
&= -2|a_2a_3|+|a_2|^2+|a_3|^2 \\ 
&= |a_2|^2-2|a_2||a_3|+|a_3|^2 = (|a_2|-|a_3|)^2 \ge 0. 
\end{align*}
So $-\det[d(Fe_3F^*)]$ is positive definite on $W$. Set $f=Fe_3F^*$.  
Then $f|_W:W\to\mathbb{S}^3_1$ determines a conformal immersion with 
induced metric 
\[
ds^2=f^*ds^2_{\mathbb{S}^3_1}=\langle df,df\rangle=-\det[d(Fe_3F^*)].
\]
Furthermore, $f$ is CMC $1$ by Theorem \ref{th:AA}.  
Also, by \eqref{eq:a1a2a3}, 
\[
-\det[d(FF^*)]=2|a_1|^2+|a_2|^2+|a_3|^2
\]
is positive definite on $M$.  Thus if we set 
\[
\beta =\frac{\det[d(FF^*)]}{\det[d(Fe_3F^*)]}
\]
on $W$, $\beta$ is a positive function on $W$ such that 
\[
\beta ds^2=-\det[d(FF^*)]
\]
extends to a Riemannian metric on $M$. Also, 
\[
\partial f\cdot f^{-1}=dF\cdot F^{-1}=F(F^{-1}dF)F^{-1}\ne 0, 
\]
and so $df\ne 0$.  This completes the proof.
\end{proof}

\begin{remark}
Even if $F$ is a holomorphic null immersion, 
$f=Fe_3F^*$ might not be a CMC $1$ face. 
For example, for the holomorphic null immersion 
\[
F:\mathbb{C}\ni z \mapsto 
  \begin{pmatrix} z+1 & -z \\ z & -z+1 \end{pmatrix} \in SL(2,\mathbb{C}),
\]
$f=Fe_3F^*$ degenerates everywhere on $\mathbb{C}$. 
Note that $\det[d(Fe_3F^*)]$ is identically zero here. 
\end{remark}

Using Propositions \ref{pr:UYmax2.3} and \ref{pr:UYmax2.4}, we can now extend 
the representation of Aiyama-Akutagawa for simply-connected CMC $1$ immersions 
to the case of CMC $1$ faces with possibly non-simply-connected domains. 

\begin{theorem} 
Let $M$ be a Riemann surface with a base point $z_0\in M$. Let $g$ be a 
meromorphic function and $\omega$ a holomorphic $1$-form on the universal 
cover $\widetilde{M}$ such that $d\hat s^2$ in Equation \eqref{eq:dshat^2} is 
a Riemannian metric on $\widetilde{M}$ and $|g|$ is not identically $1$. 
Choose the holomorphic immersion $F=(F_{jk}):\widetilde M\to SL(2,\mathbb{C})$ 
so that $F(z_0)=e_0$ and $F$ satisfies Equation \eqref{eq:F^-1dF}. 
Then $f:\widetilde M\to\mathbb{S}^3_1$ defined by 
Equation \eqref{eq:f=Fe_3F^*} is a CMC $1$ face that is conformal away from 
its singularities. 
The induced metric $ds^2$ on $M$, the second fundamental form $h$, 
and the hyperbolic Gauss map $G$ of $f$ are given as in 
Equation \eqref{eq:ds^2-h-G}. 
The singularities of the CMC $1$ face occur at points where $|g|=1$. 

Conversely, 
let $M$ be a Riemann surface and $f:M\to\mathbb{S}^3_1$ a CMC $1$ face. 
Then there exists a meromorphic function $g$ 
$($so that $|g|$ is not identically $1)$ and holomorphic $1$-form $\omega$ on 
$\widetilde{M}$ such that $d\hat s^2$ is a Riemannian metric on 
$\widetilde{M}$, and such that Equation \eqref{eq:f=Fe_3F^*} holds, 
where $F:\widetilde{M}\to SL(2,\mathbb{C})$ is an immersion which satisfies 
Equation \eqref{eq:F^-1dF}. 
\end{theorem}

\begin{proof}
First we prove the first paragraph of the theorem.  Since 
\[
d\left( Fe_3F^*\right)
=F\left(F^{-1}dF\cdot e_3+(F^{-1}dF\cdot e_3)^*\right) F^*,
\]
we have 
\[
-\det [d( Fe_3F^*)]=(1-|g|^2)^2\omega\bar\omega .
\]
Also, since $d\hat s^2$ gives a Riemannian metric on $\widetilde{M}$,  
$\omega$ has a zero of order $k$ if and only if $g$ has a pole of order 
$k/2\in\mathbb{N}$.  
Therefore $\det [d( Fe_3F^*)]=0$ if and only if $|g|=1$. 
Hence by Proposition \ref{pr:UYmax2.4}, 
$f=Fe_3F^*:\widetilde{M}\to\mathbb{S}^3_1$ is a CMC $1$ face, 
and $p\in\widetilde{M}$ is a singular point of $f$ if and only if $|g(p)|=1$, 
proving the first half of the theorem.  

We now prove the second paragraph of the theorem.  
By Proposition \ref{pr:UYmax2.3}, there exists a holomorphic null 
lift $F:\widetilde{M}\to SL(2,\mathbb{C})$ of the CMC $1$ face $f$.  
Then by the same argument as in the proof of Proposition \ref{pr:UYmax2.4}, 
we may set 
\[
F^{-1}dF=\begin{pmatrix}a_1&a_2\\a_3&-a_1\end{pmatrix}, 
\]
where $a_j$ ($j=1,2,3$) are holomorphic $1$-forms such that \eqref{eq:a1a2a3} 
holds. 
By changing $F$ into $FB$ for some constant $B\in SU(1,1)$, if necessary, 
we may assume that $a_3$ is not identically zero. 
We set 
\[
\omega :=a_3,\qquad g:=\frac{a_1}{a_3}.
\]
Then $\omega$ is a holomorphic $1$-form and $g$ is a meromorphic function. 
Since $a_2/a_3=a_2a_3/a_3^2=-(a_1/a_3)^2=-g^2$, 
we see that Equation \eqref{eq:F^-1dF} holds.  
Since $|a_2|-|a_3|$ is not identically zero 
(by the same argument as in the proof of Proposition \ref{pr:UYmax2.4}), 
$|g|$ is not identically one. 
Also, since $g^2\omega =a_2$ is holomorphic, 
\eqref{eq:d(FF^*)} implies 
$-\det [d(FF^*)]=(1+|g|^2)^2\omega\bar\omega =d\hat s^2$ is 
positive definite, so $d\hat s^2$ gives a Riemannian 
metric on $\widetilde{M}$, proving the converse part of the theorem. 
\end{proof}

\begin{remark}
Since the hyperbolic Gauss map $G$ has the geometric meaning explained in 
\eqref{rm:hypG} of Remark \ref{rm:AArep}, $G$ is single-valued on $M$ itself, 
although $F$ might not be. 
We also note that, by Equation \eqref{eq:ds^2-h-G}, the Hopf differential 
$Q$ is single-valued on $M$ as well. 
\end{remark}

\begin{remark}\label{rm:FtoFB}
Let $F$ be a holomorphic null lift of a CMC $1$ face $f$ with Weierstrass 
data $(g,\omega)$.  For any constant matrix 
\begin{equation}\label{eq:BinSU11}
B=\begin{pmatrix}\bar p&-q\\-\bar q&p\end{pmatrix}\in SU(1,1),
\qquad p\bar p-q\bar q=1,
\end{equation}
$FB$ is also a holomorphic null lift of $f$.  The Weierstrass data 
$(\hat g,\hat\omega)$ corresponding to $(FB)^{-1}d(FB)$ is given by 
\begin{equation}\label{eq:hatghatomega}
\hat g=\frac{pg+q}{\bar qg +\bar p}\quad\text{and}\quad
\hat\omega =(\bar qg+\bar p)^2\omega .
\end{equation}
Two Weierstrass data $(g,\omega)$ and $(\hat g,\hat\omega)$ are called 
{\em equivalent} if they satisfy Equation \eqref{eq:hatghatomega} for some 
$B$ as in Equation \eqref{eq:BinSU11}.  See Equation (1.6) in \cite{UY1}.  
We shall call the equivalence class of the Weierstrass data $(g,\omega)$ 
the Weierstrass data associated to $f$.  
When we wish to emphasis that $(g,\omega)$ is determined by $F$, not $FB$ for 
some $B$, we call $(g,\omega)$ the Weierstrass data associated to $F$.  

On the other hand, the Hopf differential $Q$ and the hyperbolic Gauss map $G$ 
are independent of the choice of $F$, because  
\[
\hat\omega d\hat g=\omega dg \quad\text{and}\quad
\frac{\bar p dF_{11}+\bar q dF_{12}}{\bar p dF_{21}+\bar q dF_{22}}
=\frac{dF_{11}}{dF_{21}},\quad\text{where}\quad F=(F_{jk}). 
\]
This can also be seen from \eqref{rm:hypG} of Remark \ref{rm:AArep}, 
which implies that $G$ is determined just by $f$.  Then $S(g)=S(\hat g)$ and 
\eqref{rm:2Q=Sg-SG} of Remark \ref{rm:AArep} imply $Q$ is independent of the 
choice of $F$ as well. 
\end{remark}

\section{CMC $1$ faces with elliptic ends} 
\label{sec:unitaryends} 

It is known that the only complete spacelike CMC $1$ immersion is 
a flat totally umbilic immersion \cite{Ak, R} 
(see Example \ref{ex:horos} below). 
In the case of non-immersed CMC $1$ faces in $\mathbb{S}^3_1$, we now define 
the notions of completeness and finiteness of total curvature away from 
singularities, like in \cite{KUY, UYmax}. 

\begin{definition}\label{df:comp-fin}
Let $M$ be a Riemann surface and $f:M\to\mathbb{S}^3_1$ a CMC $1$ face. 
Set $ds^2=f^*(ds^2_{\mathbb{S}^3_1})$. 
$f$ is {\em complete} (resp. of {\em finite type}) if there exists a compact 
set $C$ and a symmetric $(0,2)$-tensor $T$ on $M$ such that $T$ vanishes on 
$M\setminus C$ and $ds^2+T$ is a complete (resp. finite total curvature) 
Riemannian metric. 
\end{definition}

\begin{remark}
For CMC $1$ immersions in $\mathbb{S}^3_1$, the Gauss curvature $K$ is 
non-negative. So the total curvature is the same as the total absolute 
curvature. 
However, for CMC $1$ faces with singular points the total curvature is never 
finite, not even on neighborhoods of singular points, as can be seen from the 
form $K=4dgd\bar g/(1-|g|^2)^4\omega\bar\omega$ for the Gaussian curvature, 
see also \cite{ER}.  
Hence the phrase ``finite type'' is more appropriate in Definition 
\ref{df:comp-fin}. 
\end{remark}

\begin{remark}
The universal covering of a complete (resp. finite type) CMC $1$ face might 
not be complete (resp. finite type), because the singular set might not be 
compact on the universal cover. 
\end{remark}

Let $f:M\to\mathbb{S}^3_1$ be a complete CMC $1$ face of finite type. 
Then $(M,ds^2+T)$ is a complete Riemannian manifold of finite total curvature. 
So by \cite[Theorem 13]{H}, $M$ has finite topology, 
where we define a manifold to be of finite topology if it is diffeomorphic to 
a compact manifold with finitely many points removed.  
The {\em ends} of $f$ correspond to the removed points of that Riemann surface.

Let $\varrho :\widetilde{M}\to M$ be the universal cover of $M$, and 
$F:\widetilde{M}\to SL(2,\mathbb{C})$ a holomorphic null lift of a CMC $1$ 
face $f:M\to\mathbb{S}^3_1$.  We fix a point $z_0\in M$.  
Let $\gamma :[0,1]\to M$ be a loop so that $\gamma (0)=\gamma (1)=z_0$. 
Then there exists a unique deck transformation $\tau$ of $\widetilde{M}$ 
associated to the homotopy class of $\gamma$. 
We define the monodromy representation $\Phi_\gamma$ of $F$ as 
\[
F\circ\tau =F\Phi_\gamma .
\]
Note that $\Phi_\gamma\in SU(1,1)$ for any loop $\gamma$, since $f$ is 
well-defined on $M$. 
So $\Phi_\gamma$ is conjugate to either 
\begin{equation}\label{eq:e-h-p}
\mathcal{E}=\begin{pmatrix}e^{i\theta}&0\\0&e^{-i\theta}\end{pmatrix}
\quad\text{or}\quad
\mathcal{H}=\pm\begin{pmatrix}\cosh s&\sinh s\\\sinh s&\cosh s\end{pmatrix}
\quad\text{or}\quad
\mathcal{P}=\pm\begin{pmatrix}1+i&1\\1&1-i\end{pmatrix}
\end{equation}
for $\theta\in [0,2\pi)$, $s\in\mathbb{R}\setminus\{0\}$.  

\begin{definition}\label{df:e-h-p}
Let $f:M\to\mathbb{S}^3_1$ be a complete CMC $1$ face of finite type with 
holomorphic null lift $F$. 
An end of $f$ is called an {\em elliptic end} or {\em hyperbolic end} or 
{\em parabolic end} if its monodromy representation is 
conjugate to $\mathcal{E}$ or $\mathcal{H}$ or $\mathcal{P}$ 
in $SU(1,1)$, respectively.  
\end{definition}

\begin{remark}
A matrix 
\[
X=\begin{pmatrix}p&q\\\bar q&\bar p\end{pmatrix}\in SU(1,1)
\]
acts on the hyperbolic plane in the Poincar\'e model 
$\mathbb{H}^2=(\{w\in\mathbb{C}\,|\,|w|<1\},
               ds^2_{\mathbb{H}^2}=4dwd\bar w/(1-|w|^2)^2)$ 
as an isometry: 
\[ 
\mathbb{H}^2\ni w\mapsto\frac{pw+q}{\bar qw+\bar p}\in\mathbb{H}^2.
\]
$X$ is called {\em elliptic} if this action has only one fixed point 
which is in $\mathbb{H}^2$. 
$X$ is called {\em hyperbolic} if there exist two fixed points, 
both in the ideal boundary $\partial\mathbb{H}^2$. 
$X$ is called {\em parabolic} if there exists only one fixed point 
which is in $\partial\mathbb{H}^2$. 
This is what motivates the terminology in Definition \ref{df:e-h-p}. 
\end{remark}

Since any matrix in $SU(2)$ is conjugate to $\mathcal{E}$ in $SU(2)$, 
CMC $1$ immersions in $\mathbb{H}^3$ and CMC $1$ faces with elliptic ends in 
$\mathbb{S}^3_1$ share many analogous properties.  
So in this paper we consider CMC $1$ faces with only elliptic ends. 
We leave the study of hyperbolic ends and parabolic ends for another 
occasion. 

\begin{proposition}\label{pr:dshat2}
Let $V$ be a neighborhood of an end of $f$ and $f|_V$ a spacelike 
CMC $1$ immersion of finite total curvature which is complete at the end.  
Suppose the end is elliptic. 
Then there exists 
a holomorphic null lift $F:\widetilde{V}\to SL(2,\mathbb{C})$ of $f$ 
with Weierstrass data $(g,\omega)$ associated to $F$ 
such that 
\[
d\hat s^2|_V=(1+|g|^2)^2\omega\bar\omega
\]
is single-valued on $V$. 
Moreover, $d\hat s^2|_V$ has finite total curvature and is complete at the 
end. 
\end{proposition}

\begin{proof} 
Let $\gamma :[0,1]\to V$ be a loop around the end and $\tau$ the 
deck transformation associated to $\gamma$. 
Take a holomorphic null lift $F_0:\widetilde{V}\to SL(2,\mathbb{C})$ of $f$. 
Then by definition of an elliptic end, there exists a $\theta\in [0,2\pi )$ 
such that 
\[
F_0\circ\tau =F_0PE_\theta P^{-1},
\]
where $P\in SU(1,1)$ and 
\[
E_\theta =\begin{pmatrix}e^{i\theta}&0\\0&e^{-i\theta}\end{pmatrix}.
\]
Defining the holomorphic null lift $F=F_0P$ of $f$ and defining 
$(g,\omega)$ to be the Weierstrass data associated to $F$, 
we have 
\[
F\circ\tau =FE_\theta ,\quad
g\circ\tau =e^{-2i\theta}g\quad\text{and}\quad
\omega\circ\tau =e^{2i\theta}\omega .
\]
Thus, $|g\circ\tau|=|g|$ and $|\omega\circ\tau|=|\omega|$. 
This implies $d\hat s^2|_V$ is single-valued on $V$. 
Let $T$ be a $(0,2)$-tensor as in Definition \ref{df:comp-fin}. 
Then by Equation \eqref{eq:ds^2-h-G}, 
we have $ds^2+T\le d\hat s^2|_V$ on $V\setminus C$. 
Thus, if $ds^2+T$ is complete, $d\hat s^2|_V$ is also complete. 
We denote the Gaussian curvature of the metric $d\hat s^2|_V$ by 
$K_{d\hat s^2|_V}$ (note that $K_{d\hat s^2|_V}$ is non-positive).  
Then we have 
\[
(-K_{d\hat s^2|_V})d\hat s^2|_V
=\frac{4dgd\bar g}{(1+|g|^2)^2}
\le\frac{4dgd\bar g}{(1-|g|^2)^2}=Kds^2
\]
on $V\setminus C$. 
Thus, if $ds^2+T$ is of finite total curvature, the total curvature of 
$d\hat s^2|_V$ is finite, proving the proposition. 
\end{proof}

\begin{proposition}\label{co:UYmax4.5}
Let $f:M\to\mathbb{S}^3_1$ be a complete CMC $1$ face of finite type with 
elliptic ends. 
Then there exists a compact Riemann surface $\overline{M}$ and finite number 
of points $p_1,\dots ,p_n\in\overline{M}$ so that $M$ is biholomorphic to 
$\overline{M}\setminus\{p_1,\dots ,p_n\}$. 
Moreover, the Hopf differential $Q$ of $f$ extends meromorphically to 
$\overline{M}$. 
\end{proposition}

\begin{proof}
Since $f$ is of finite type, $M$ is finitely connected, 
by \cite[Theorem 13]{H}. 
Consequently, there exists a compact region $M_0\subset M$, bounded by a 
finite number of regular Jordan curves $\gamma_1,\dots,\gamma_n$, such that 
each component $M_j$ of $M\setminus M_0$ can be conformally mapped onto the 
annulus $D_j=\{z\in\mathbb{C}\,|\,r_j<|z|<1\}$, where $\gamma_j$ corresponds 
to the set $\{|z|=1\}$. 
Then by Proposition \ref{pr:dshat2}, there exists $d\hat s^2|_{M_j}$ which is 
single-valued on $M_j$ and is of finite total curvature and is complete at 
the end, and so that $K_{d\hat s^2|_{M_j}}$ is non-positive.  
Therefore by \cite[Proposition III-16]{L} or \cite[Theorem 9.1]{Os}, 
$r_j=0$, and hence each $M_j$ is biholomorphic to the punctured disk 
$\{z\in\mathbb{C}\,|\,0<|z|<1\}$. 
We can, using the biholomorphism from $M_j$ to $D_j$, replace $M_j$ in $M$ 
with $D_j$ without affecting the conformal structure of $M$.  
Thus, without loss of generality, $M=\overline{M}\setminus\{p_1,\dots,p_n\}$ 
for some compact Riemann surface $\overline{M}$ and a finite number of points 
$p_1,\dots,p_n$ in $\overline{M}$, and each $M_j$ becomes a punctured disk 
about $p_j$.  
Hence, by \eqref{rm:FF*} of Remark \ref{rm:AArep}, we can apply 
\cite[Proposition 5]{B} to $\hat f_j:=\hat f|_{M_j}$ to see that $Q=\hat Q$ 
extends meromorphically to $M_j\cup\{p_j\}$, proving the proposition. 
\end{proof}

Let $\Delta =\{z\in\mathbb{C}\,|\,|z|<1\}$ and 
$\Delta^*=\Delta\setminus\{0\}$.  Let $f:\Delta^*\to\mathbb{S}^3_1$ be a 
conformal spacelike CMC $1$ immersion of finite total curvature which has a 
complete elliptic end at the origin.  
Then we can take the Weierstrass data associated to $f$ in the following form: 
\begin{equation}\label{eq:UY1(W1)}
\left\{\begin{array}{lll}
       g=z^\mu\tilde g(z),          & \mu >0,     & \tilde g(0)\ne 0, \\
  \omega=w(z)dz=z^\nu\tilde w(z)dz, & \nu \le -1, & \tilde w(0)\ne 0,
       \end{array}\right.
\end{equation}
where $\tilde g$ and $\tilde w$ are holomorphic functions on $\Delta$ and 
$\mu +\nu\in\mathbb{Z}$. 
(See \cite{UY1} and also \cite[Proposition 4]{B} for case that 
$\mu, \nu\in\mathbb{R}$. 
Then applying a transformation as in Equation \eqref{eq:hatghatomega} if 
necessary, we may assume $\mu >0$. 
Completeness of the end then gives $\nu\le -1$.) 

\begin{definition}
The end $z=0$ of $f:\Delta^*\to\mathbb{S}^3_1$ is called {\em regular} if the 
hyperbolic Gauss map $G$ extends meromorphically across the end. 
Otherwise, the end is called {\em irregular}. 
\end{definition}

Since $Q$ extends meromorphically to each end, 
we have the following proposition, 
by \eqref{rm:FF*} and \eqref{rm:2Q=Sg-SG} of Remark \ref{rm:AArep}: 

\begin{proposition}\label{pr:regularity} {\rm \cite[Proposition 6]{B}}
An end $f:\Delta^*\to\mathbb{S}^3_1$ is regular if and only if the order at 
the end of the Hopf differential of $f$ is at least $-2$.
\end{proposition}

\begin{remark}
In \cite{LY}, Lee and Yang define {\em normal ends} and {\em abnormal ends}.  
Both normal and abnormal ends are biholomorphic to a punctured disk 
$\Delta^*$, and the Hopf differential has a pole of order $2$ at the origin. 
Normal ends are elliptic ends, and abnormal ends are hyperbolic ends. 
Moreover, the Lee-Yang catenoids with normal ends are complete in our sense.  
However, the Lee-Yang catenoids with abnormal ends include incomplete 
examples, because the singular set of these examples accumulates at the ends. 
(In fact, CMC $1$ face with hyperbolic ends cannot be complete, 
 which is shown in \cite{FRUYY}).  
\end{remark}

\section{Embeddedness of elliptic ends} 
\label{sec:embeddedness} 

In this section we give a criterion for embeddedness of elliptic ends, 
which is based on results in \cite{UY1}. 
This criterion will be applied in the next section. 

Let $f:\Delta^*\to\mathbb{S}^3_1$ be a conformal spacelike CMC $1$ immersion 
of finite total curvature with a complete regular elliptic end 
at the origin. 
Let $\gamma :[0,1]\to\Delta^*$ be a loop around the origin and 
$\tau$ the deck transformation of $\widetilde{\Delta^*}$ 
associated to the homotopy class of $\gamma$. 
Then by the same argument as in the proof of Proposition \ref{pr:dshat2}, 
there exists the holomorphic null lift $F$ of $f$ such that 
$F\circ\tau =FE_\theta$ for some $\theta\in [0,2\pi)$, where 
$E_\theta =\mathrm{diag}(e^{i\theta},e^{-i\theta})$. 
Since $E_\theta\in SU(2)$, 
$\hat f:=FF^*$ in $\mathbb{H}^3$ is single-valued on $\Delta^*$. 

Let $(g,\omega)$ be the Weierstrass data associated to $F$, 
defined as in \eqref{eq:UY1(W1)}. 
Then by \eqref{rm:FF*} in Remark \ref{rm:AArep}, 
$\hat f:\Delta^*\to\mathbb{H}^3$ is a conformal CMC $1$ immersion with the 
induced metric $d\hat s^2=(1+|g|^2)^2\omega\bar\omega$. 
Thus by the final sentence of Proposition \ref{pr:dshat2}, 
$\hat f$ has finite total curvature and is complete at the origin. 

Since $\hat f$ has the same Hopf differential $Q$ as $f$, 
$f$ having a regular end immediately implies that $\hat f$ has a regular end. 

Furthermore we show the following theorem:

\begin{theorem}\label{lm:sub2}
Let $f:\Delta^*\to\mathbb{S}^3_1$ be a conformal spacelike CMC $1$ immersion 
of finite total curvature with a complete regular elliptic end 
at the origin. 
Then there exists a holomorphic null lift 
$F:\widetilde{\Delta^*}\to SL(2,\mathbb{C})$ of $f$ $($that is, $f=Fe_3F^*)$ 
such that $\hat f=FF^*$ is a conformal CMC $1$ finite-total-curvature 
immersion from $\Delta^*$ into $\mathbb{H}^3$ with a complete regular end 
at the origin.  
Moreover, $f$ has an embedded end if and only if $\hat f$ has an embedded end. 
\end{theorem}

\begin{remark}
The converse of the first part of Theorem \ref{lm:sub2} is also true, 
that is, the following holds: 
Let $\hat f:\Delta^*\to\mathbb{H}^3$ be a conformal CMC $1$ immersion 
of finite total curvature with a complete regular end at the origin. 
Take a holomorphic null lift $F:\widetilde{\Delta^*}\to SL(2,\mathbb{C})$ of 
$\hat f$ (that is, $f=FF^*$) such that the associated Weierstrass data 
$(g,\omega)$ is written as in \eqref{eq:UY1(W1)}.  
Then $f=Fe_3F^*$ is a conformal spacelike CMC $1$ finite-total-curvature 
immersion from $\Delta^*$ into $\mathbb{S}^3_1$ with a complete regular 
elliptic end at the origin.  
Moreover, $\hat f$ has an embedded end if and only if $f$ has an embedded end. 
See Proposition \ref{pr:UYmax5.4} below. 
\end{remark}

We already know that such an $F$ exists. 
So we must prove the equivalency of embeddedness between the ends of 
$f$ and $\hat f$. 
To prove this we prepare three lemmas. 

\begin{lemma}[{\cite[Lemma 5.3]{UY1}}]\label{lm:UY1Lemma5.3} 
There exists a $\Lambda\in SL(2,\mathbb{C})$ such that 
\begin{equation}
\Lambda F=\begin{pmatrix}z^{\lambda_1}a(z)&z^{\lambda_2}b(z)\\
                     z^{\lambda_1-m_1}c(z)&z^{\lambda_2-m_2}d(z)\end{pmatrix},
\label{eq:LambdaF}
\end{equation}
where $a$, $b$, $c$, $d$ are holomorphic functions on $\Delta$ that do not 
vanish at the origin, and $\lambda_j\in\mathbb{R}$ and $m_j\in\mathbb{N}$ 
$(j=1,2)$ are defined as follows:
\begin{enumerate}
\item If $\mathrm{Ord}_0Q=\mu +\nu-1=-2$, then 
\begin{equation}\label{eq:UY1(5.8)}
m_1=m_2, \quad 
\lambda_1=\frac{-\mu +m_j}{2}\quad\text{and}\quad
\lambda_2=\frac{\mu +m_j}{2}. 
\end{equation}
\item If $\mathrm{Ord}_0Q=\mu +\nu-1\ge -1$, then 
\begin{equation}\label{eq:UY1(5.9)}
m_1=-(\nu +1), \quad m_2=2\mu +\nu +1, \quad 
\lambda_1=0\quad\text{and}\quad\lambda_2=m_2. 
\end{equation}
\end{enumerate}
\end{lemma}

Note that in either case we have $\nu <-1$ and 
\begin{equation}\label{eq:lambda-m}
\lambda_1<\lambda_2, \quad 
\lambda_1-m_1<\lambda_2-m_2\quad\text{and}\quad\lambda_1-m_1<0. 
\end{equation}
Note also that in the second case we have $m_1<m_2$. 

\begin{proof}[Proof of Lemma \ref{lm:UY1Lemma5.3}]
$F$ satisfies Equation \eqref{eq:F^-1dF}, which is precisely Equation (1.5) 
in \cite{UY1}. 
So we can apply \cite[Lemma 5.3]{UY1}, since that lemma is based 
on Equation (1.5) in \cite{UY1}.  This gives the result.  
\end{proof}

It follows that 
\begin{align*}
\Lambda f\Lambda^*
&= (\Lambda F)e_3(\Lambda F)^* \\
&= \begin{pmatrix}
   |z|^{2\lambda_1}|a|^2-|z|^{2\lambda_2}|b|^2 &
   |z|^{2\lambda_1}\bar z^{-m_1}a\bar c-|z|^{2\lambda_2}\bar z^{-m_2}b\bar d \\
   |z|^{2\lambda_1}z^{-m_1}\bar ac-|z|^{2\lambda_2}z^{-m_2}\bar bd &
   |z|^{2(\lambda_1-m_1)}|c|^2-|z|^{2(\lambda_2-m_2)}|d|^2
   \end{pmatrix}. 
\end{align*}
Note that $\Lambda f\Lambda^*$ is congruent to $f=Fe_3F^*$, 
because $(\Lambda F)^{-1}d(\Lambda F)=F^{-1}dF$ determines both the 
first and second fundamental forms as in Equation \eqref{eq:ds^2-h-G}.  

\begin{figure}[htbp] 
\begin{center}
\unitlength=1.0pt
\begin{picture}(300,200)
\put(158,196){\makebox(0,0)[cc]{$x_0$}}
\put(270,92){\makebox(0,0)[rc]{$(x_1,x_2,x_3)$}}
\put(148,50){\makebox(0,0)[rc]{$(-1,0,0,0)$}}
\put(148,150){\makebox(0,0)[rc]{$(1,0,0,0)$}}
\put(191,95){\makebox(0,0)[ct]
 {$\underbrace{\phantom{M}}_{\mathbb{D}^3}$}}
\put(109,95){\makebox(0,0)[ct]
 {$\underbrace{\phantom{M}}_{\mathbb{D}^3}$}}
\put(250,20){\makebox(0,0)[cc]{$\mathbb{S}^3_1$}}
\put(150,50){\circle*{2}}
\put(150,150){\circle*{2}}
\put(182,100){\circle*{2}}
\put(118,100){\circle*{2}}
\put(50,100){\vector(1,0){200}}
\put(150,0){\vector(0,1){200}}
\bezier70(150,100)(100,150)(50,200)
\bezier70(150,100)(200,150)(250,200)
\bezier400(50,190)(150,100)(50,10)
\bezier400(250,190)(150,100)(250,10)
\bezier300(150,50)(182,99)(214,148)
\bezier34(150,150)(183,150)(216,150)
\end{picture}
\end{center}
\caption{The projection $p$.}%
\label{fg:proj}
\end{figure}
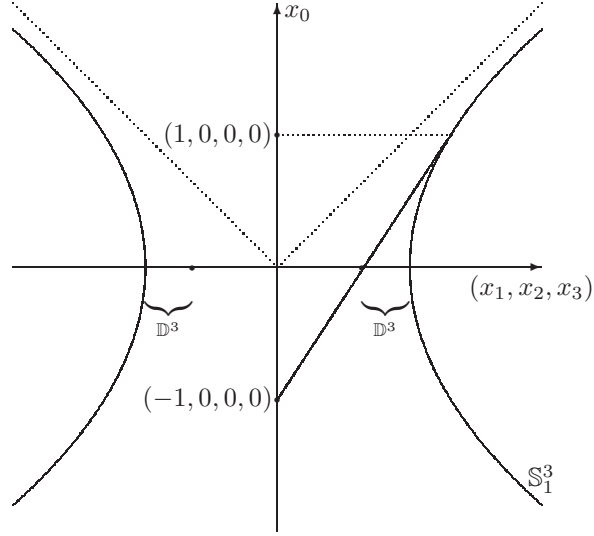

To study the behavior of the elliptic end $f$, we present the elliptic end in 
a $3$-ball model as follows: 
We set 
\[
\begin{pmatrix}x_0+ x_3&x_1+ix_2\\
               x_1-ix_2&x_0- x_3\end{pmatrix}
=\Lambda f\Lambda^*. 
\]
Since 
\begin{align*}
x_0
&= \frac{1}{2}\text{trace}(\Lambda f\Lambda^*) \\
&= \frac{1}{2}
   \left(|z|^{2\lambda_1}|a|^2-|z|^{2\lambda_2}|b|^2
        +|z|^{2(\lambda_1-m_1)}|c|^2-|z|^{2(\lambda_2-m_2)}|d|^2\right), 
\end{align*}
\eqref{eq:lambda-m} implies that $\lim_{z\to 0}x_0(z)=\infty$.  
So we may assume that $x_0(z)>1$ for any $z\in\Delta^*$. 
So we can define a bijective map 
\[
p:\{(x_0,x_1,x_2,x_3)\in\mathbb{S}^3_1\,|\,x_0>1\}
  \to\mathbb{D}^3:=D^3_1\setminus \overline{D^3_{1/\sqrt{2}}}
\]
as 
\[
p(x_0,x_1,x_2,x_3):=\frac{1}{1+x_0}(x_1,x_2,x_3),
\]
where $D^3_r$ denotes the open ball of radius $r$ in $\mathbb{R}^3$ 
and $\overline{D^3_r}=D^3_r\cup \partial D^3_r$ 
(See Figure \ref{fg:proj}). 

We set $(X_1,X_2,X_3)=p\circ (\Lambda f\Lambda^*)$.  
Then we have 
\begin{align}
X_1+iX_2
&= \frac{2a\bar c|z|^{2(\lambda_1-m_1)}z^{m_1}
         -2b\bar d|z|^{2(\lambda_2-m_2)}z^{m_2}}
         {2+|a|^2|z|^{2\lambda_1}+|c|^2|z|^{2(\lambda_1-m_1)}
           -|b|^2|z|^{2\lambda_2}-|d|^2|z|^{2(\lambda_2-m_2)}}, 
\label{eq:X1+iX2}\\
X_3
&= \frac{|a|^2|z|^{2\lambda_1}-|c|^2|z|^{2(\lambda_1-m_1)}
         -|b|^2|z|^{2\lambda_2}+|d|^2|z|^{2(\lambda_2-m_2)}}
         {2+|a|^2|z|^{2\lambda_1}+|c|^2|z|^{2(\lambda_1-m_1)}
           -|b|^2|z|^{2\lambda_2}-|d|^2|z|^{2(\lambda_2-m_2)}}. 
\label{eq:X3}
\end{align}

We now define a function $U:\Delta^*\to\mathbb{C}$ that will be useful for 
proving Theorem \ref{lm:sub2}: 
\begin{equation}\label{eq:UY1(5.15)}
U(z)=z^{-m_1}(X_1+iX_2). 
\end{equation}
Then by making just a few sign changes to the argument 
in \cite[Lemma 5.4]{UY1}, we have the following lemma: 

\begin{lemma}\label{lm:UY1Lemma5.4}
The function $U$ is a $C^\infty$ function on $\Delta^*$ that extends 
continuously to $\Delta$ with $U(0)\ne 0$. Moreover, 
\begin{equation}\label{eq:zdU->0}
\lim_{z\to 0}z\frac{\partial U}{\partial z}=0\quad\text{and}\quad
\lim_{z\to 0}z\frac{\partial U}{\partial\bar z}=0. 
\end{equation}
\end{lemma}

\begin{proof}
If $\mathrm{Ord}_0Q=-2$, then Equation \eqref{eq:UY1(5.15)} is reduced to 
\begin{equation}\label{eq:UY1(5.19)}
U(z)=\frac{2a\bar c-2b\bar d|z|^{2\mu}}
          {2|z|^{\mu+m_1}+|a|^2|z|^{2m_1}
           +|c|^2-|b|^2|z|^{2(\mu +m_1)}-|d|^2|z|^{2\mu}}.
\end{equation}
Then we have $U(0)=2a(0)/c(0)\ne 0$, because $\mu>0$ and $m_1\ge 1$. 
Also, by a straightforward calculation, 
we see that Equation \eqref{eq:zdU->0} holds. 

If $\mathrm{Ord}_0Q\ge -1$, then Equation \eqref{eq:UY1(5.15)} is reduced to 
\begin{equation}\label{eq:UY1(5.20)}
U(z)=\frac{2a\bar c-2b\bar dz^{m_2}\bar z^{m_1}}
          {(2+|a|^2-|d|^2)|z|^{2m_1}-|b|^2|z|^{2(m_1+m_2)}+|c|^2}.
\end{equation}
Then we have $U(0)=2a(0)/c(0)\ne 0$, because $m_1,m_2\ge 1$. 
Also, since $U$ is $C^1$-differentiable at the origin, 
we see that Equation \eqref{eq:zdU->0} holds. 
\end{proof}

Also, we have the following lemma, 
analogous to Remark 5.5 in \cite{UY1}: 

\begin{lemma}\label{lm:UY1Remark5.5}
$p\circ (\Lambda f\Lambda^*)=(X_1,X_2,X_3)$ converges 
to the single point $(0,0,-1)\in\partial D^3_1$. 
Moreover, $p\circ (\Lambda f\Lambda^*)$ is tangent to $\partial D^3_1$ 
at the end $z=0$. 
\end{lemma}

\begin{proof}
By \eqref{eq:UY1(5.8)}--\eqref{eq:X3}, we see that 
\[
\lim_{z\to 0}(X_1,X_2,X_3)=(0,0,-1). 
\]
Define a function $V:\Delta^*\to\mathbb{C}$ by $V:=z^{-m_1}(X_3+1)$. 
Then from either case of the proof of Lemma \ref{lm:UY1Lemma5.4}, we see that 
\[
\lim_{z\to 0}V=0\quad\text{and}\quad
\lim_{z\to 0}z\frac{\partial V}{\partial z}=0. 
\]
Therefore 
\begin{align*}
0&= \lim_{z\to 0}z\frac{\partial V}{\partial z}
  = \lim_{z\to 0}\left(z^{-m_1+1}\frac{\partial X_3}{\partial z}\right)
                -m_1V(0), \\ 
0&= \lim_{z\to 0}z\frac{\partial U}{\partial z}
  = \lim_{z\to 0}\left(z^{-m_1+1}\frac{\partial (X_1+iX_2)}{\partial z}\right)
                -m_1U(0) 
\end{align*}
imply that 
\[
\lim_{X_1+iX_2\to 0}\frac{\partial X_3}{\partial (X_1+iX_2)}
=\frac{V(0)}{U(0)}=0, 
\]
proving the lemma. 
\end{proof}

\begin{proof}[Proof of Theorem \ref{lm:sub2}]
\cite[Theorem 5.2]{UY1} shows that $\hat f$ has an embedded end 
if and only if $m_1=1$, so the theorem will be proven by showing that 
$f$ also has an embedded end if and only if $m_1=1$.  We now show this:  

Since $U(0)\ne 0$, by taking a suitable branch we can define the function 
$u:\Delta\to\mathbb{C}$ by 
\[
u(z)=z\left(U(z)\right)^{1/m_1}. 
\]
Then $u$ is a $C^1$ function such that 
\begin{equation}\label{eq:UY1(5.21)}
\frac{\partial u}{\partial z}(0)\ne 0\quad\text{and}\quad
\frac{\partial u}{\partial\bar z}=0. 
\end{equation}

Assume that $m_1=1$.  By Equation \eqref{eq:UY1(5.15)}, $X_1+iX_2=u$ holds. 
Then by Equation \eqref{eq:UY1(5.21)}, 
$X_1+iX_2$ is one-to-one on some neighborhood of the origin $z=0$. 
Hence $f$ has an embedded end. 

Conversely, assume that $f$ has an embedded end. 
Let $p_3:\mathbb{D}^3\to\Delta$ be the projection defined 
by $p_3(X_1,X_2,X_3)=X_1+iX_2$. 
By Equation \eqref{eq:UY1(5.15)} we have 
\begin{equation}\label{eq:UY1(5.22)}
p_3\circ p\circ (\Lambda f\Lambda^*)=u^{m_1}. 
\end{equation}
By Equations \eqref{eq:UY1(5.21)} and \eqref{eq:UY1(5.22)} the map 
$p_3\circ p\circ (\Lambda f\Lambda^*)$ is an $m_1$-fold cover of 
$\Delta^*_\varepsilon =\{z\in\mathbb{C}\,|\,0<|z|<\varepsilon\}$ for a 
sufficiently small $\varepsilon >0$. 
Thus, by Lemma \ref{lm:UY1Remark5.5}, 
$m_1$ must be $1$, by the same topological arguments as at the end of the 
proof of Theorem 5.2 in \cite{UY1}. 

Therefore, we have that $f$ has an embedded end if and only if $m_1=1$. 
\end{proof}

\section{The Osserman-type inequality} 
\label{sec:Ossermanineq} 

Let $f:M=\overline{M}\setminus\{p_1,\dots ,p_n\}\to\mathbb{S}^3_1$ be a 
complete CMC $1$ face of finite type with hyperbolic Gauss map $G$ and 
Hopf differential $Q$. 

\begin{definition}
We set 
\begin{equation}
\label{eq:norm-metric}
 d\hat s^{\sharp 2}
:=(1+|G|^2)^2\frac{Q}{dG}\overline{\left(\frac{Q}{dG}\right)}
\end{equation}
and call it the {\em lift-metric} of the CMC $1$ face $f$. 
We also set 
\[
d\hat\sigma^{\sharp 2}:=(-K_{d\hat s^{\sharp 2}})d\hat s^{\sharp 2}
                       =\frac{4dGd\overline{G}}{(1+|G|^2)^2}.
\]
\end{definition}

\begin{remark}
Since $G$ and $Q$ are defined on $M$, both $d\hat s^{\sharp 2}$ and 
$d\hat\sigma^{\sharp 2}$ are also defined on $M$. 
\end{remark}

We define the order of pseudometrics as in \cite[Definition 2.1]{UY3}, 
that is: 

\begin{definition}
A pseudometric $d\varsigma^2$ on $\overline{M}$ is of order $m_j$ at $p_j$ if 
$d\varsigma^2$ has a local expression 
\[
d\varsigma^2=e^{2u_j}dzd\bar z
\]
around $p_j$ such that $u_j-m_j\log |z-z(p_j)|$ is continuous at $p_j$. 
We denote $m_j$ by $\mathrm{Ord}_{p_j}(d\varsigma^2)$. 
In particular, if $d\varsigma^2$ gives a Riemannian metric around $p_j$, 
then $\mathrm{Ord}_{p_j}(d\varsigma^2)=0$. 
\end{definition}

We now apply \cite[Lemma 3]{UY2} for regular ends in $\mathbb{H}^3$ to regular 
elliptic ends in $\mathbb{S}^3_1$, that is, we show the following proposition: 

\begin{proposition}\label{pr:UY2Lemma3}
Let $f:\Delta^*\to\mathbb{S}^3_1$ be a conformal spacelike CMC $1$ immersion 
of finite total curvature with a complete regular elliptic end 
at the origin $z=0$. 
Then the following inequality holds:
\begin{equation}\label{eq:UY2(20)}
\mathrm{Ord}_0(d\hat\sigma^{\sharp 2})-\mathrm{Ord}_0(Q)\ge 2.
\end{equation}
Moreover, equality holds if and only if the end is embedded.
\end{proposition}

\begin{proof} 
By Theorem \ref{lm:sub2}, there exists a holomorphic null lift 
$F:\widetilde{\Delta^*}\to SL(2,\mathbb{C})$ of $f$ such that  
$\hat f=FF^*:\Delta^*\to\mathbb{H}^3$ is a conformal CMC $1$ immersion  
of finite total curvature with a complete regular end at the origin.  
Then by \cite[Lemma 3]{UY2}, we have \eqref{eq:UY2(20)}. 
Moreover, equality holds if and only if the end of $\hat f$ is embedded, 
by \cite[Lemma 3]{UY2}. 
This is equivalent to the end of $f$ being embedded, 
by Theorem \ref{lm:sub2}, proving the proposition. 
\end{proof}

The following lemma is a variant on known results in \cite{Yu, KTUY}. 
In fact, \cite{Yu} showed that $d\hat s^2$ is complete if and only if 
$d\hat s^{\sharp 2}$ is complete, see also \cite[Lemma 4.1]{KTUY}. 

\begin{lemma}\label{lm:UY4.3}
Let $f:M\to\mathbb{S}^3_1$ be a CMC $1$ face. Assume that each end of $f$ is 
regular and elliptic.  
If $f$ is complete and of finite type, then the lift-metric 
$d\hat s^{\sharp 2}$ is complete and of finite total curvature on $M$. 
In particular, 
\begin{equation}
\mathrm{Ord}_{p_j}(d\hat s^{\sharp 2})\le -2 
\label{eq:dshatsharp}
\end{equation}
holds at each end $p_j$ $(j=1,\dots ,n)$. 
\end{lemma}

\begin{proof}
Since $f$ is complete and of finite type, each end is complete and has 
finite total curvature. 
Then by \eqref{eq:UY2(20)} and  the relation 
\begin{equation}\label{eq:SharpGauss-eqn}
 \mathrm{Ord}_{p_j}(d\hat s^{\sharp 2})
+\mathrm{Ord}_{p_j}(d\hat\sigma^{\sharp 2})
=\mathrm{Ord}_{p_j}(Q)
\end{equation}
(which follows from the Gauss equation 
$d\hat\sigma^{\sharp 2}d\hat s^{\sharp 2}=4Q\overline{Q}$), 
we have \eqref{eq:dshatsharp} at each end $p_j$. 
Hence $d\hat s^{\sharp 2}$ is a complete metric. 
Also, again by \eqref{eq:UY2(20)}, we have 
\[
\mathrm{Ord}_{p_j}(d\hat\sigma^{\sharp 2})\ge 2+\mathrm{Ord}_{p_j}(Q)\ge 0,
\]
because $p_j$ is regular (that is, $\mathrm{Ord}_{p_j}(Q)\ge -2$, 
by Proposition \ref{pr:regularity}). 
This implies that the total curvature of $d\hat s^{\sharp 2}$ is finite. 
\end{proof}

\begin{remark}
Consider a CMC $1$ face with regular elliptic ends. 
If it is complete and of finite type, 
then, by Lemma \ref{lm:UY4.3}, 
the lift-metric is complete and of finite total curvature. 
But the converse is not true.  See \cite{FRUYY}. 
\end{remark}

\begin{theorem}[Osserman-type inequality]\label{th:main}
Let $f:M\to\mathbb{S}^3_1$ be a complete CMC $1$ face of finite type with $n$ 
elliptic ends and no other ends.  
Let $G$ be its hyperbolic Gauss map.  Then the following inequality holds:
\begin{equation}\label{eq:Oss-ineq}
2\deg (G)\geq -\chi (M)+n, 
\end{equation}
where $\deg (G)$ is the mapping degree of $G$ $($if $G$ has essential 
singularities, then we define $\deg (G)=\infty)$. 
Furthermore, equality holds if and only if each end is regular and embedded.%
\end{theorem}

\begin{remark}
As we remarked in the introduction, the completeness of a CMC $1$ face $f$ 
implies that $f$ must be of finite type, and each end must be elliptic. 
See the forthcoming paper \cite{FRUYY}. 
\end{remark}

\begin{proof}[Proof of Theorem \ref{th:main}]
Recall that we can set $M=\overline{M}\setminus\{p_1,\ldots ,p_n\}$, 
where $\overline{M}$ is a compact Riemann surface and $p_1,\ldots ,p_n$ is a 
set of points in $\overline{M}$, by Proposition \ref{co:UYmax4.5}. 
If $f$ has irregular ends, then $G$ has essential singularities at those ends.
So $\deg (G)=\infty$ and then \eqref{eq:Oss-ineq} automatically holds. 
Hence we may assume $f$ has only regular ends.
Using the Riemann-Hurwicz formula and the Gauss equation 
$d\hat s^{\sharp 2}d\hat\sigma^{\sharp 2}=4Q\overline{Q}$, 
we have
\begin{align*}
2\deg (G)
&= \chi(\overline{M})
  +\sum_{p\in\overline{M}}\mathrm{Ord}_{p}d\hat\sigma^{\sharp 2} \\
&= \chi(\overline{M})
  +\sum_{p\in\overline{M}}
   \left(\mathrm{Ord}_{p}Q-\mathrm{Ord}_{p}d\hat s^{\sharp 2}\right) \\
&= \chi(\overline{M})
  +\sum_{p\in\overline{M}}\mathrm{Ord}_{p}Q
  -\sum_{p\in M}\mathrm{Ord}_{p}d\hat s^{\sharp 2}
  -\sum_{j=1}^n\mathrm{Ord}_{p_j}d\hat s^{\sharp 2} \\
&= -\chi(\overline{M})-\sum_{j=1}^n\mathrm{Ord}_{p_j}d\hat s^{\sharp 2} \\
&\ge -\chi(\overline{M})+2n \quad
      \text{(because $d\hat s^{\sharp 2}$ is complete, 
             by \eqref{eq:dshatsharp})} \\
&= -\chi(M)+n.
\end{align*}
Equality in \eqref{eq:Oss-ineq} holds if and only if equality 
in \eqref{eq:dshatsharp} holds at each end, which is equivalent to 
equality in \eqref{eq:UY2(20)} holding at each end, 
by Equation \eqref{eq:SharpGauss-eqn}. 
Thus by Proposition \ref{pr:UY2Lemma3}, we have the conclusion. 
\end{proof}

\section{Examples} 
\label{sec:examples} 

To visualize CMC $1$ faces, we use the {\em hollow ball model} of 
$\mathbb{S}^3_1$, as in \cite{KY, LY, Y}. 
For any point
\[
\begin{pmatrix}x_0+x_3&x_1+ix_2\\x_1-ix_2&x_0-x_3\end{pmatrix}\leftrightarrow
(x_0,x_1,x_2,x_3)\in\mathbb{S}^3_1,
\]
define
\[
y_k=\frac{e^{\arctan x_0}}{\sqrt{1+x_0^2}}x_k,
\qquad k=1,2,3.
\]
Then $e^{-\pi}<y_1^2+y_2^2+y_3^2<e^\pi$. 
The identification $(x_0,x_1,x_2,x_3)\leftrightarrow (y_1,y_2,y_3)$ is then a 
bijection from $\mathbb{S}^3_1$ to the hollow ball 
\[
\mathscr{H}=\{(y_1,y_2,y_3)\in\mathbb{R}^3\,|\,
              e^{-\pi}<y_1^2+y_2^2+y_3^2<e^\pi\}. 
\]
So $\mathbb{S}^3_1$ is identified with the hollow ball $\mathscr{H}$, and we 
show the graphics here using this identification to $\mathscr{H}$. 

We shall first introduce four basic examples, using the same Weierstrass data 
as for CMC $1$ immersions in $\mathbb{H}^3$.  

\begin{example}\label{ex:horos}
The CMC $1$ face associated to horosphere in $\mathbb{H}^3$ is given by the 
Weierstrass data $(g,\omega)=(c_1,c_2dz)$ with $c_1\in\mathbb{C}$, 
$c_2\in\mathbb{C}\setminus\{0\}$ on the Riemann surface $\mathbb{C}$. 
This CMC $1$ face has no singularities. So this example is indeed a complete 
spacelike CMC $1$ immersion. 
\cite{Ak} and \cite[Theorem 7]{R} independently showed that the only complete 
spacelike CMC $1$ immersion in $\mathbb{S}^3_1$ 
is a flat totally umbilic immersion.  
\end{example}

\begin{figure}[htbp] 
\begin{center}
\begin{tabular}{cc}
 \includegraphics[width=.30\linewidth]{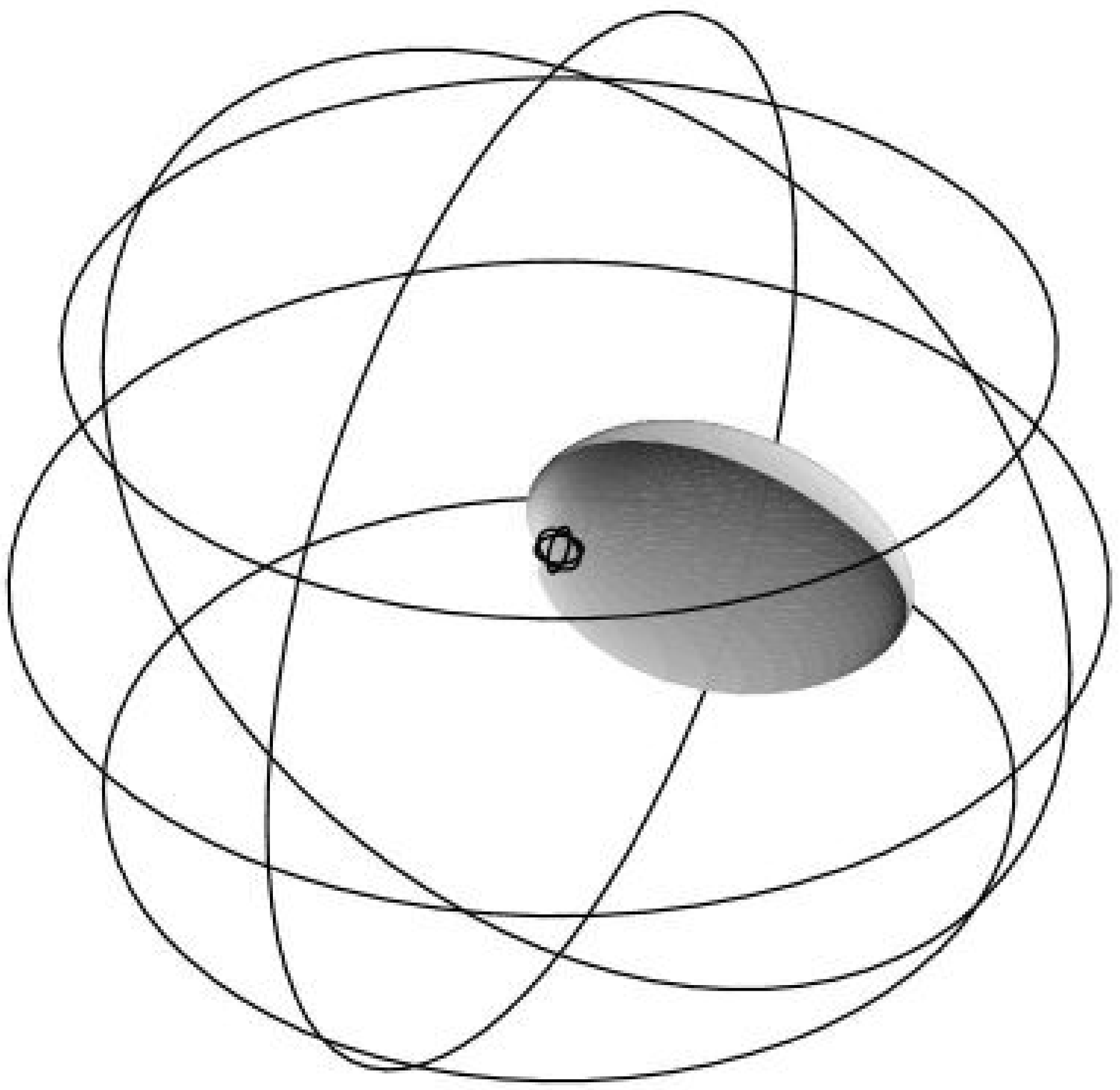} &
 \includegraphics[width=.30\linewidth]{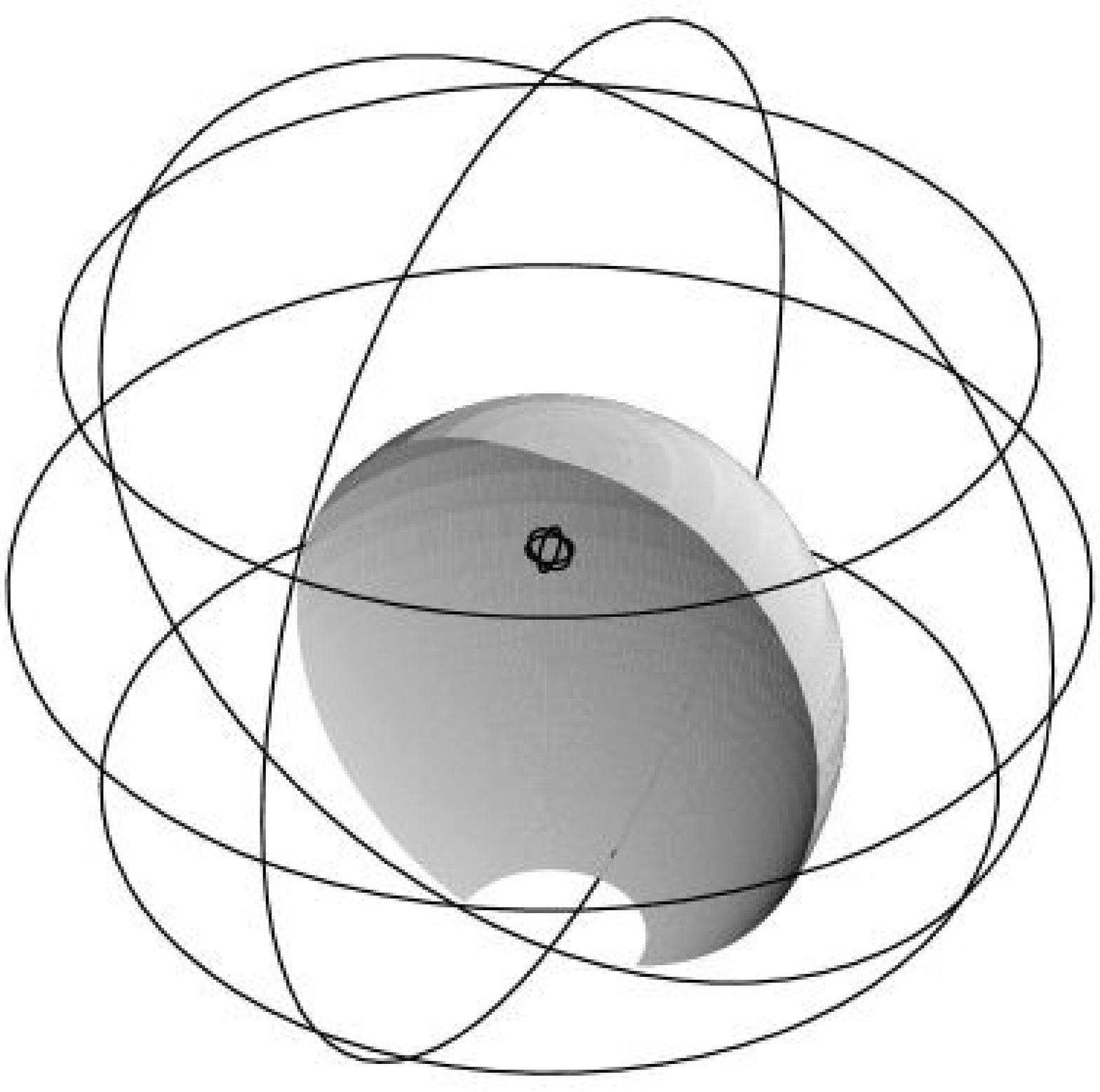} \\
 $\{z\in\mathbb{C}\,;\,
    \begin{subarray}{c}|z|<5,\\0\le\arg z\le\pi\end{subarray}\}.$ &
 $\{z\in\mathbb{C}\,;\,
    \begin{subarray}{c}|z|<10,\\0\le\arg z\le\pi\end{subarray}\}.$ 
\end{tabular}
\end{center}
\caption{Pictures of Example \ref{ex:horos}. 
         The left-hand side is drawn with $c_1=1.2$, $c_2=1$ and
         the right-hand side is drawn with $c_1=0$, $c_2=1$.}
\end{figure} 

\begin{example}\label{ex:enn}
The CMC $1$ face associated to the Enneper cousin in $\mathbb{H}^3$ is given 
by the Weierstrass data $(g,\omega)=(z,cdz)$ with 
$c\in\mathbb{R}\setminus\{0\}$ on the Riemann surface $\mathbb{C}$. 
The induced metric $ds^2=c^2(1-|z|^2)^2dzd\bar z$ degenerates where $|z|=1$. 
Take a $p\in\mathbb{C}$ which satisfies $|p|=1$. Define 
\[
\beta :=\left(\frac{1+|z|^2}{1-|z|^2}\right)^2.
\]
Then 
\[
\lim_{\begin{subarray}{c}z\to p\\z\in W\end{subarray}}\beta ds^2
=4c^2dzd\bar z. 
\]
So all singularities are admissible and hence this is a CMC $1$ face.  
Moreover, it is easily seen that this CMC $1$ face is complete and of 
finite type.  
Since this CMC $1$ face is simply-connected, the end of this CMC $1$ face is 
an elliptic end.  
Since $\mathrm{Ord}_\infty Q=-4<-2$, the end of this CMC $1$ face is 
irregular. 
Hence this CMC $1$ face does not satisfy equality in the 
inequality \eqref{eq:Oss-ineq}. 
\end{example}

\begin{figure}[htbp] 
\begin{center}
\begin{tabular}{cc}
 \includegraphics[width=.40\linewidth]{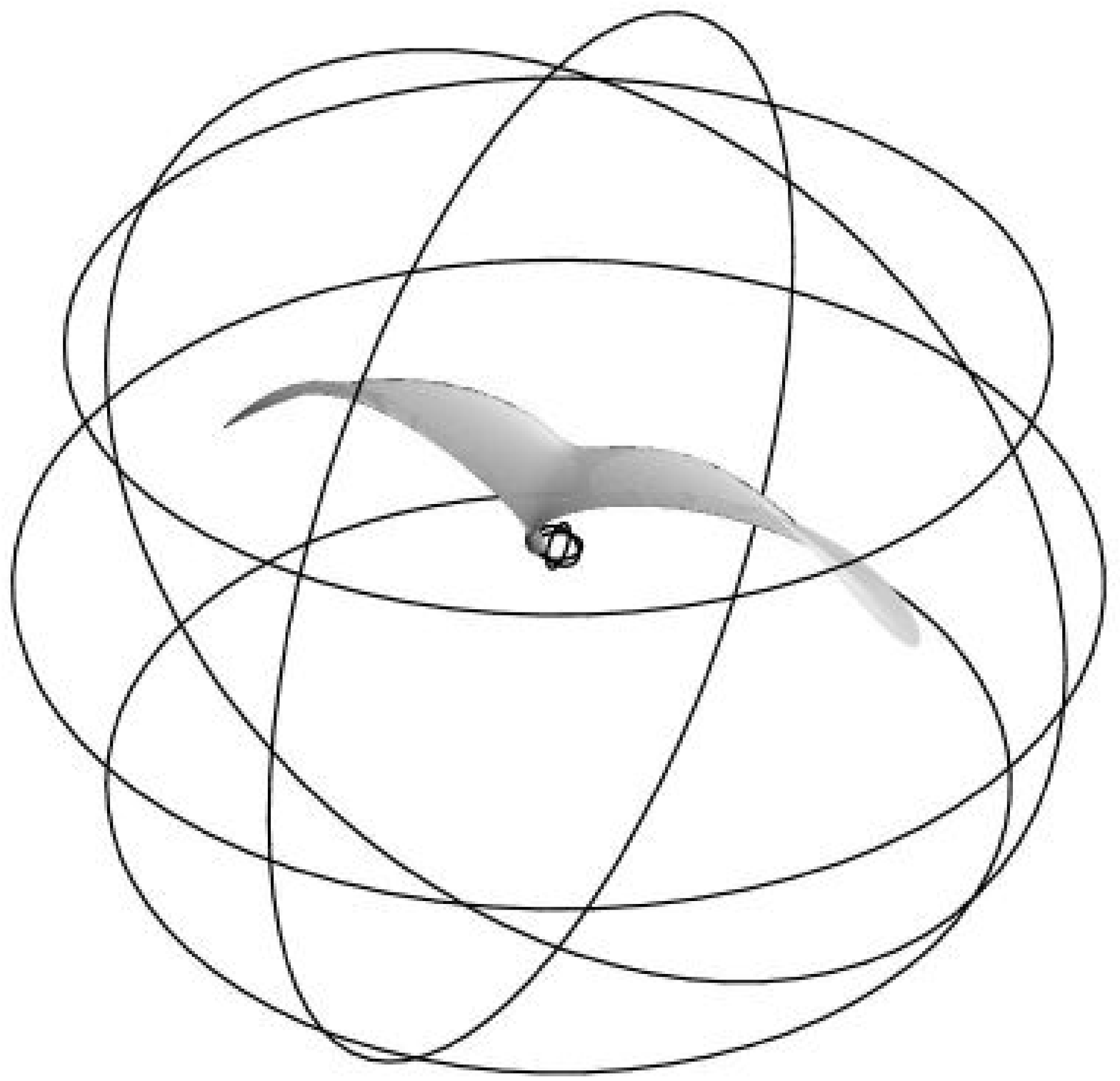} &
 \includegraphics[width=.35\linewidth]{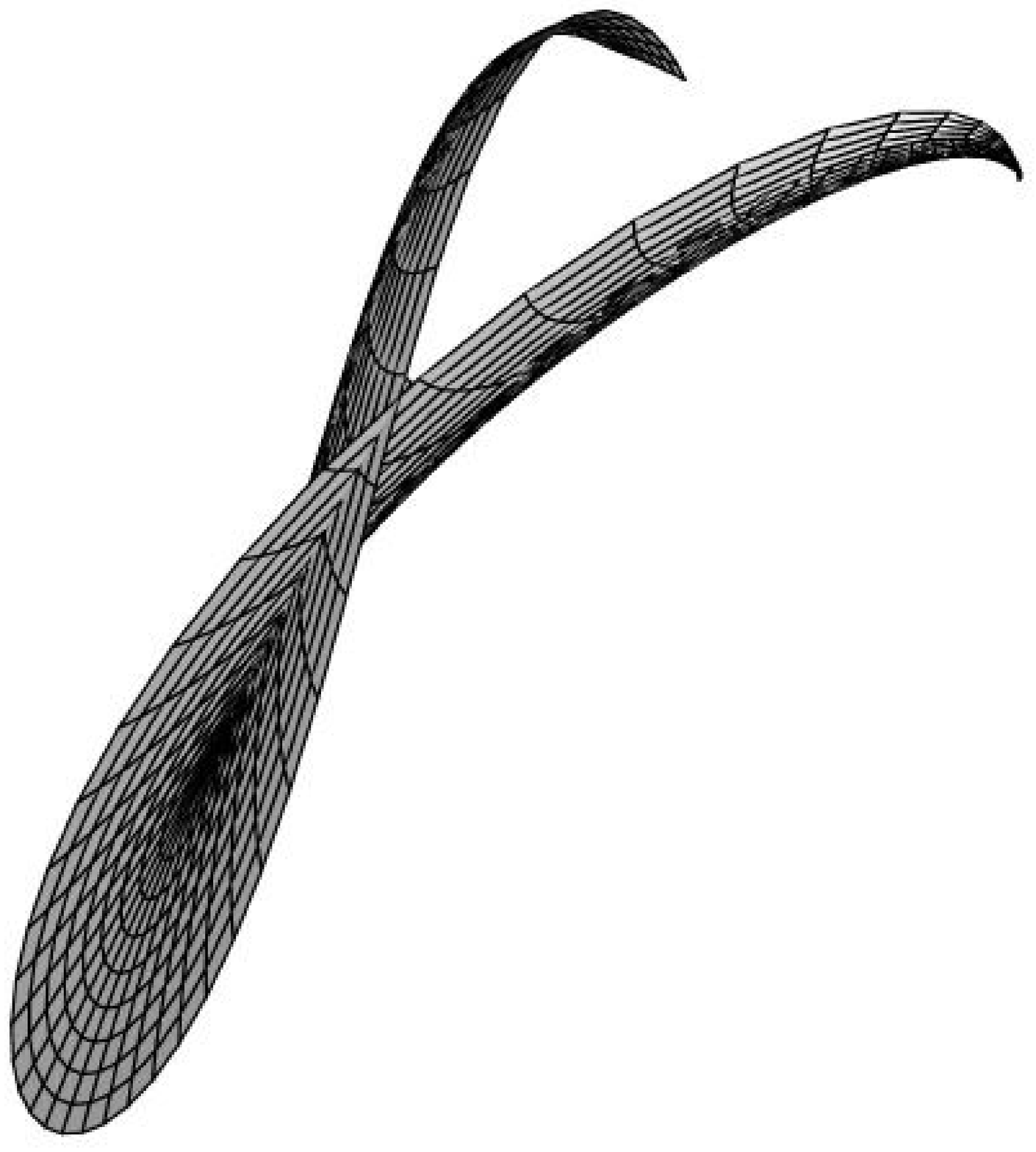} \\
 $\{z\in\mathbb{C}\,;\,|z|<1.3\}.$ &
 $\{z\in\mathbb{C}\,;\,
    \begin{subarray}{c}0.8<|z|<1.3\\\pi -1<\arg z<\pi +1\end{subarray}\}.$
\end{tabular}
\end{center}
\caption{Pictures of Example \ref{ex:enn}, where $c=1$.}
\end{figure} 

\begin{example}\label{ex:hel}
The CMC $1$ face associated to the helicoid cousin in $\mathbb{H}^3$ is given 
by the Weierstrass data $(g,\omega)=(e^z,ice^{-z}dz)$ with 
$c\in\mathbb{R}\setminus\{0\}$ on the Riemann surface $\mathbb{C}$. 
Set $z=x+iy$. The induced metric $ds^2=4c^2\sinh^2x(dx^2+dy^2)$ degenerates 
where $x=0$. 
Take a $p\in\mathbb{C}$ which satisfies $\mathrm{Re}(p)=0$. 
Define $\beta :=\tanh^{-2}x$.  Then 
\[
 \lim_{\begin{subarray}{c}z\to p\\z\in W\end{subarray}}\beta ds^2
=4c^2(dx^2+dy^2). 
\]
So all singularities are admissible and hence this is a CMC $1$ face.  
Since the singular set is non-compact, this CMC $1$ face is neither complete 
nor of finite type. 
\end{example}

\begin{figure}[htbp] 
\begin{center}
\begin{tabular}{cc}
 \includegraphics[width=.40\linewidth]{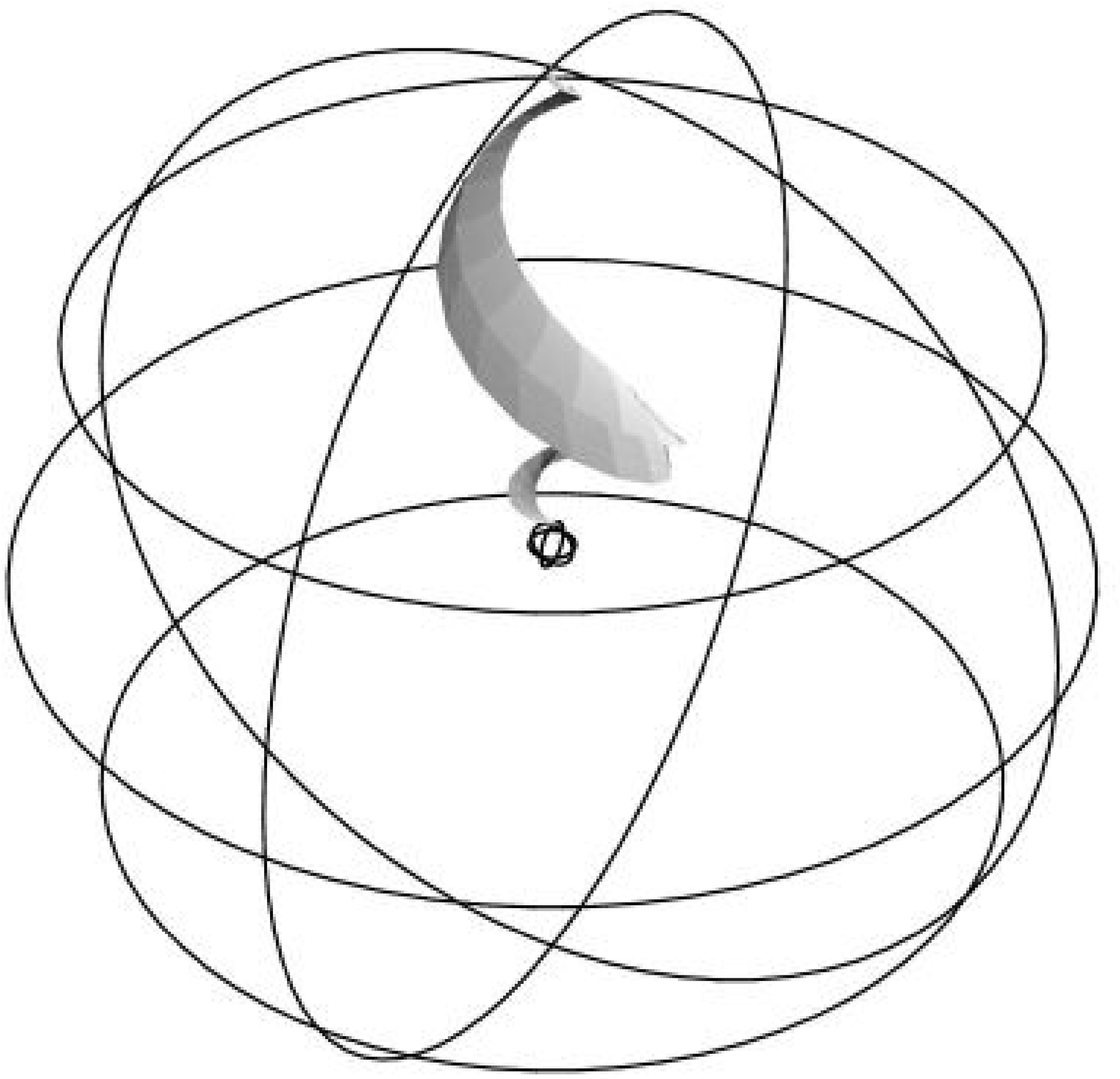} &
 \raisebox{35pt}{\includegraphics[width=.40\linewidth]{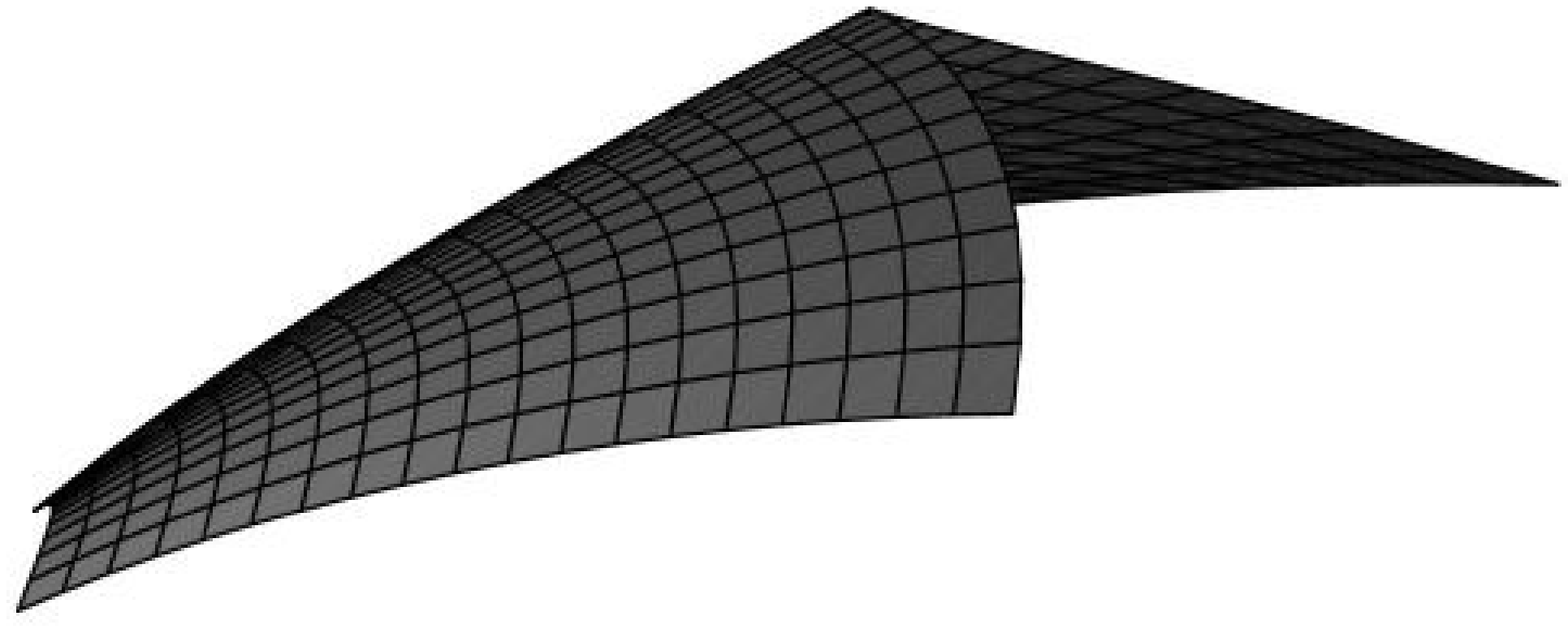}} \\
 $\{z\in\mathbb{C}\,;\,
    \begin{subarray}{c}-0.9<\mathrm{Re}z<0.9\\
                       -4\pi <\mathrm{Im}z<4\pi\end{subarray}\}.$ &
 $\{z\in\mathbb{C}\,;\,
    \begin{subarray}{c}-0.8<\mathrm{Re}z<0.8\\
                       -0.3<\mathrm{Im}z<0.3\end{subarray}\}.$
\end{tabular}
\end{center}
\caption{Pictures of Example \ref{ex:hel}, where $c=1$.}
\end{figure} 

\begin{example}\label{ex:cat}
The CMC $1$ face associated to the catenoid cousin in $\mathbb{H}^3$ is given 
by the Weierstrass data $(g,\omega)=(z^\mu,(1-\mu^2)dz/4\mu z^{\mu+1})$ with  
$\mu\in\mathbb{R}^+\setminus\{1\}$ on the Riemann surface 
$\mathbb{C}\setminus\{0\}$.  The induced metric 
\[
ds^2=
\left(\frac{|z|^{\mu}-|z|^{-\mu}}{|z|}\cdot\frac{1-\mu^2}{4\mu}\right)^2
dzd\bar z
\]
degenerates where $|z|=1$. 
Take a $p\in\mathbb{C}$ which satisfies $|p|=1$. Define 
\[
\beta :=\left(\frac{|z|^{\mu}+|z|^{-\mu}}{|z|^{\mu}-|z|^{-\mu}}\right)^2.
\]
Then 
\[
\lim_{\begin{subarray}{c}z\to p\\z\in W\end{subarray}}\beta ds^2
=4\left(\frac{1-\mu^2}{4\mu}\right)^2dzd\bar z.
\]
So all singularities are admissible and hence this is a CMC $1$ face.  
Moreover, it is easily seen that this CMC $1$ face is complete and of 
finite type.  
Since the eigenvalues of the monodromy representation at each end are 
$-e^{\mu\pi i}$, $-e^{-\mu\pi i}\in\mathbb{S}^1$, each end of this CMC $1$ 
face is an elliptic end.  
Since $\mathrm{Ord}_0 Q=\mathrm{Ord}_\infty Q=-2$, each end of this CMC $1$ 
face is regular. 
Also, each end is of this CMC $1$ face is embedded, so this CMC $1$ face 
satisfies equality in the inequality \eqref{eq:Oss-ineq}.
\end{example}

\begin{figure}[htbp] 
\begin{center}
\begin{tabular}{ccc}
 \includegraphics[width=.30\linewidth]{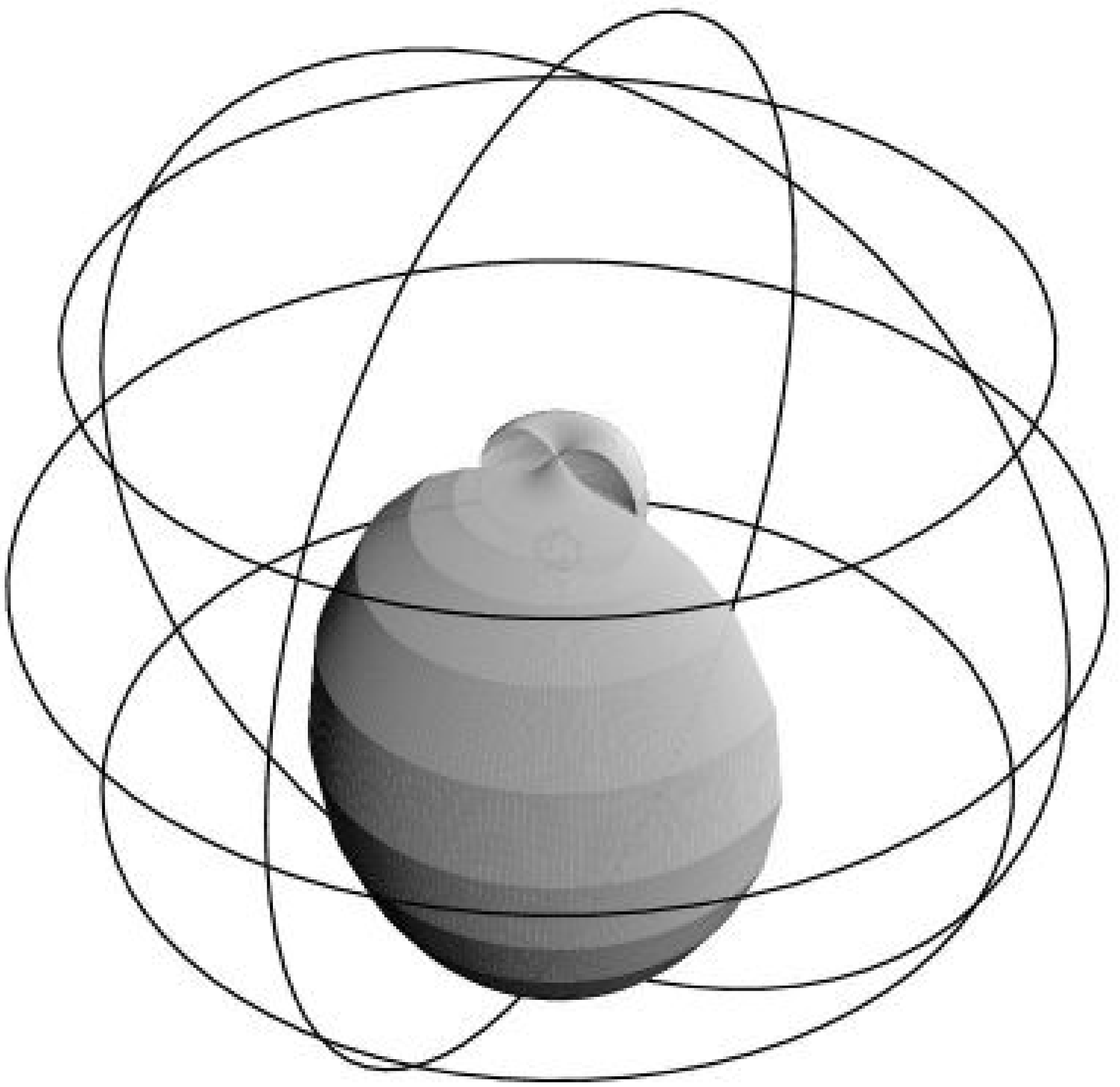} &
 \includegraphics[width=.30\linewidth]{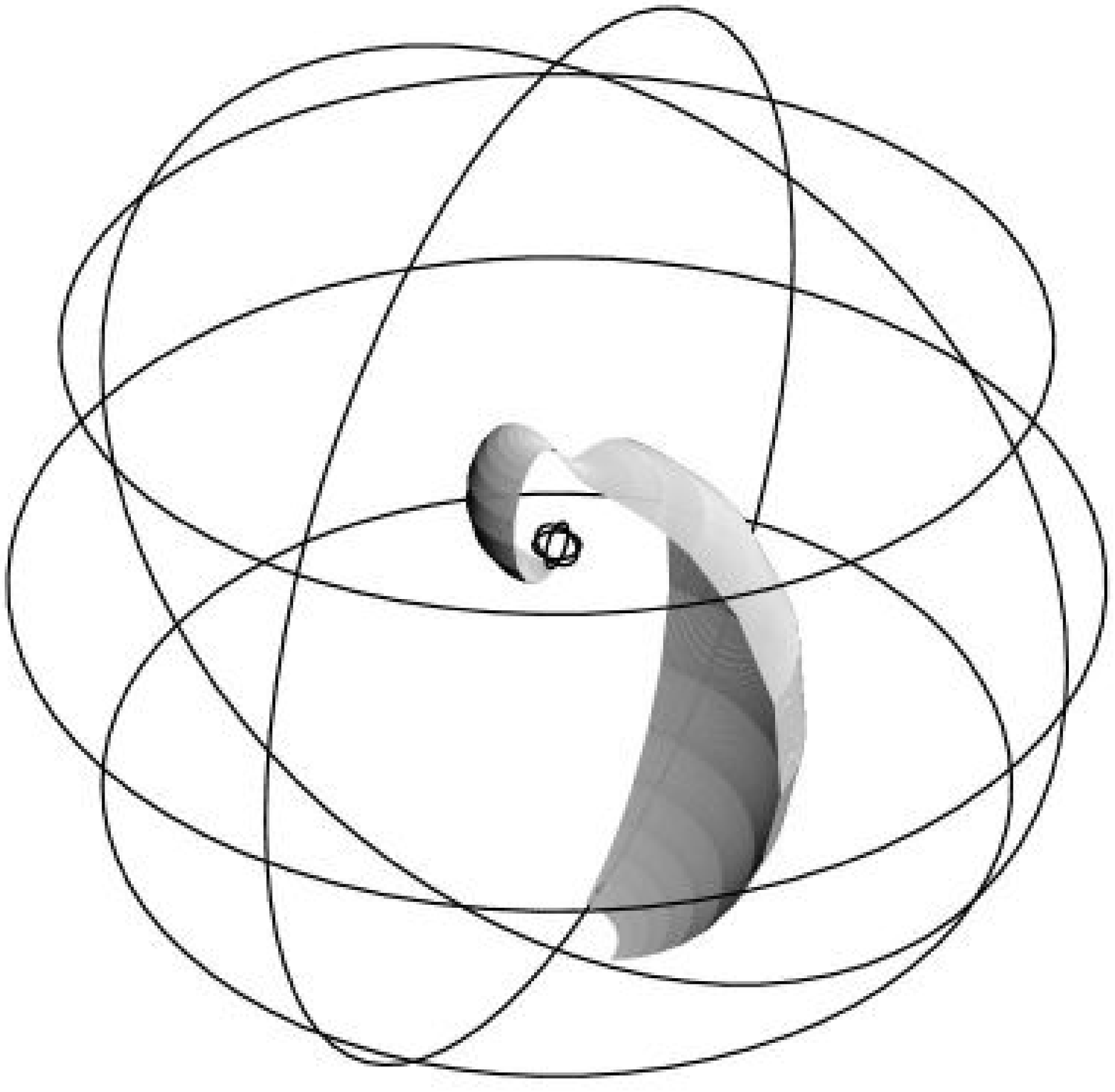} &
 \includegraphics[width=.30\linewidth]{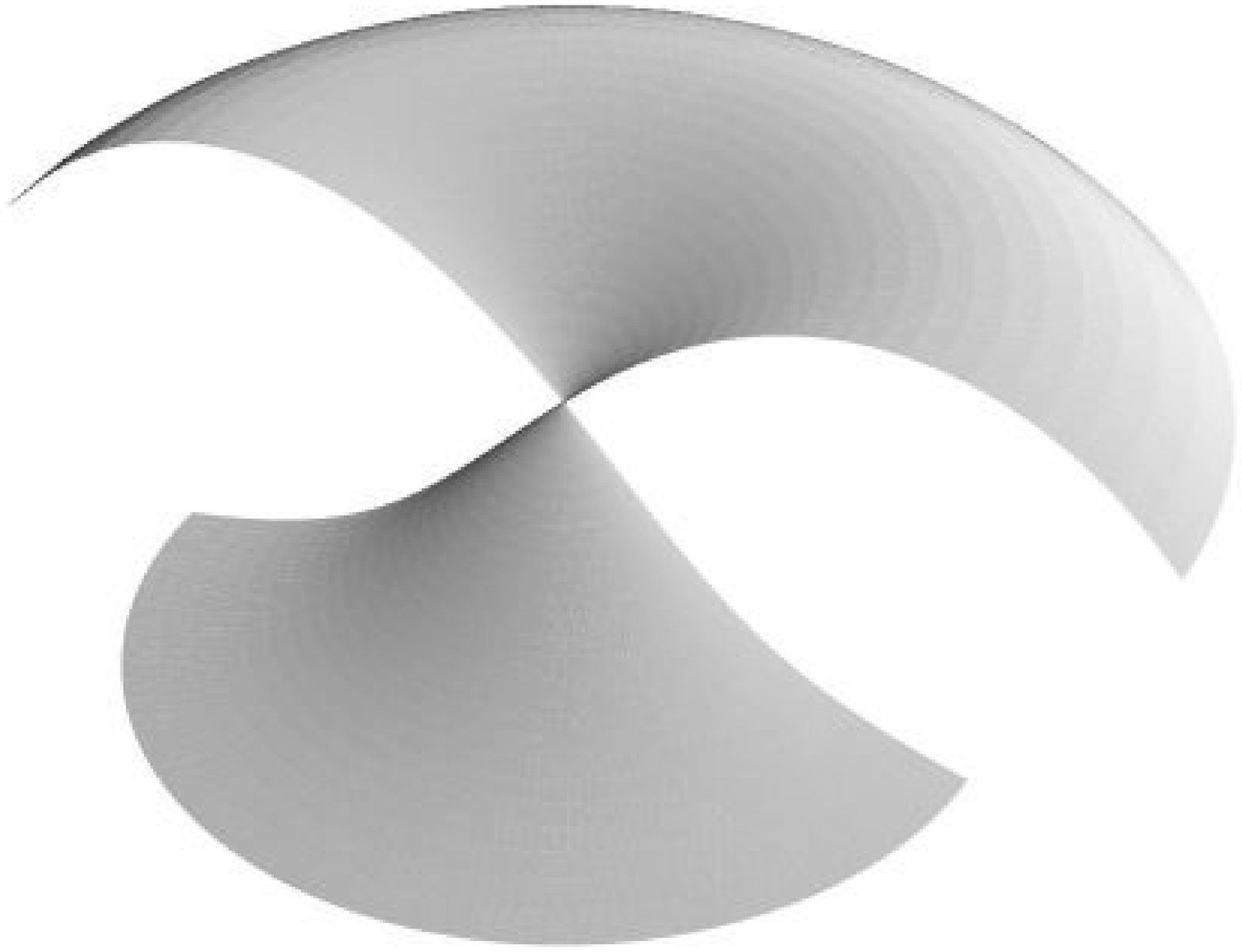} \\
 $\{z\in\mathbb{C}\,;\,
    \begin{subarray}{c}e^{-5}<|z|<e^5\\0<\arg z<\pi\end{subarray}\}.$ &
 $\{z\in\mathbb{C}\,;\,
    \begin{subarray}{c}e^{-5}<|z|<e^5\\\pi <\arg z<(3/2)\pi\end{subarray}\}.$ &
 $\{z\in\mathbb{C}\,;\,
    \begin{subarray}{c}e^{-2}<|z|<e^2\\0<\arg z<\pi\end{subarray}\}.$
\end{tabular}
\end{center}
\caption{Pictures of Example \ref{ex:cat}, where $\mu =0.8$.}
\end{figure} 

\begin{figure}[htbp] 
\begin{center}
\begin{tabular}{ccc}
 \includegraphics[width=.30\linewidth]{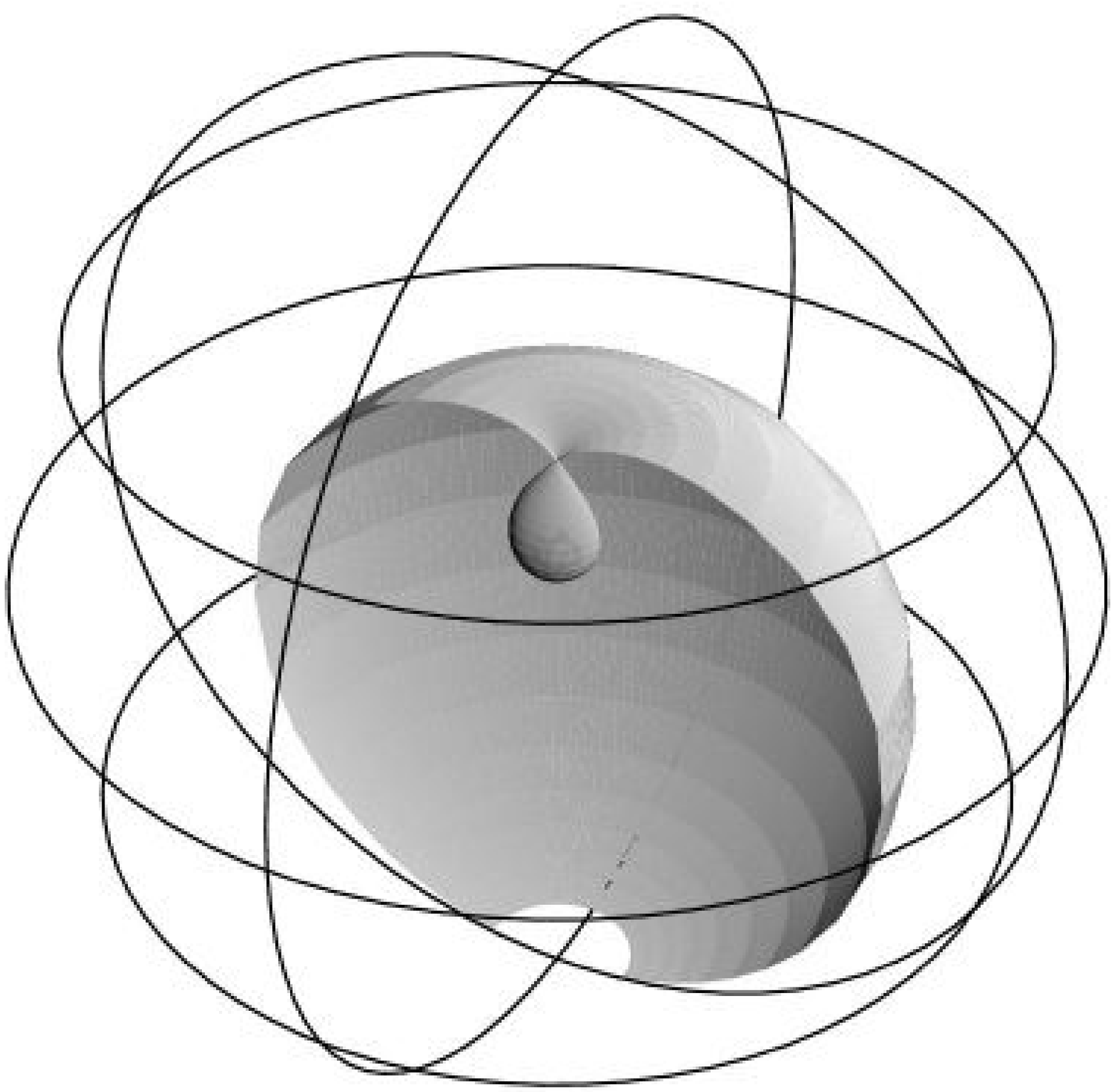} &
 \includegraphics[width=.30\linewidth]{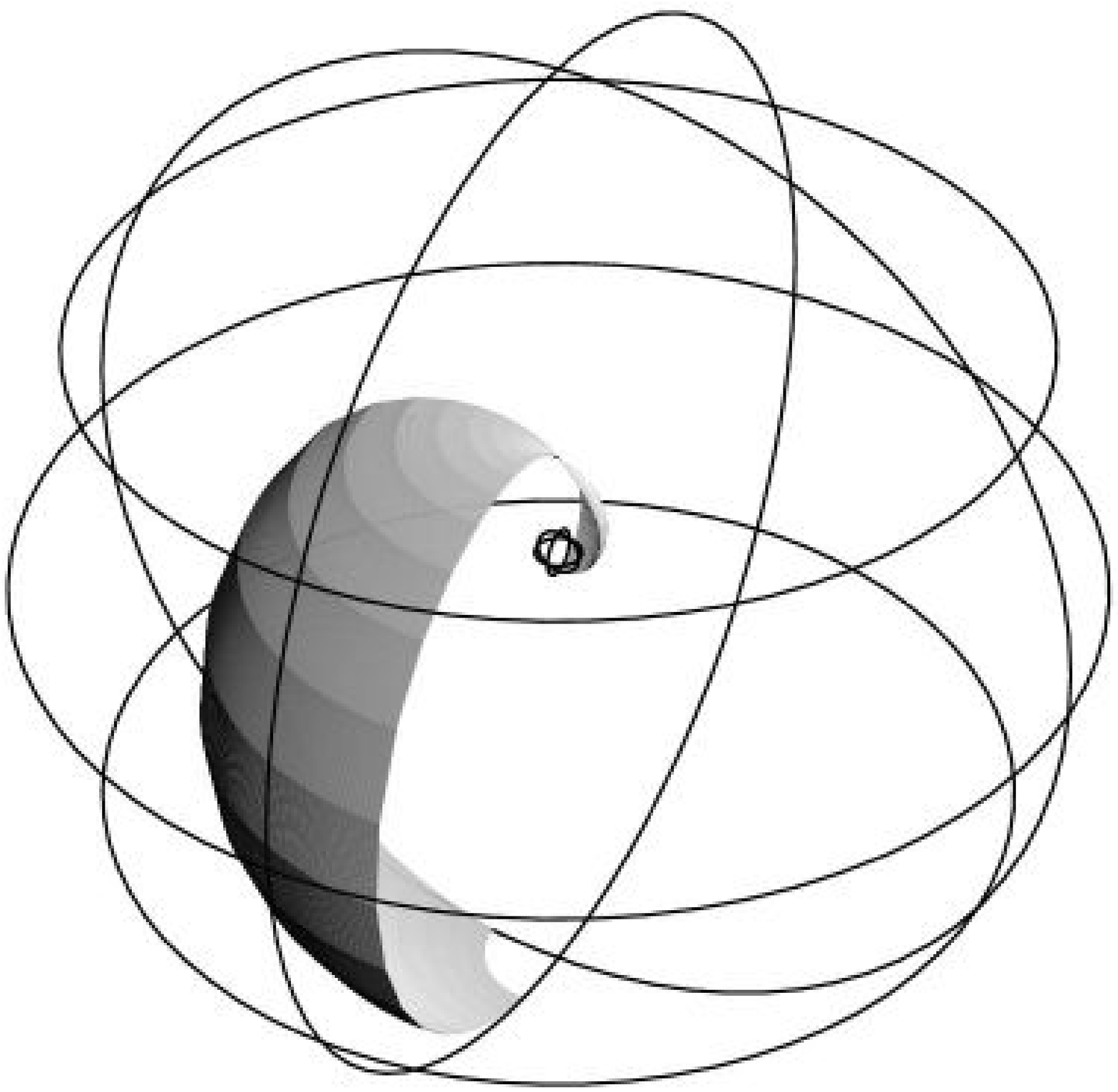} &
 \includegraphics[width=.30\linewidth]{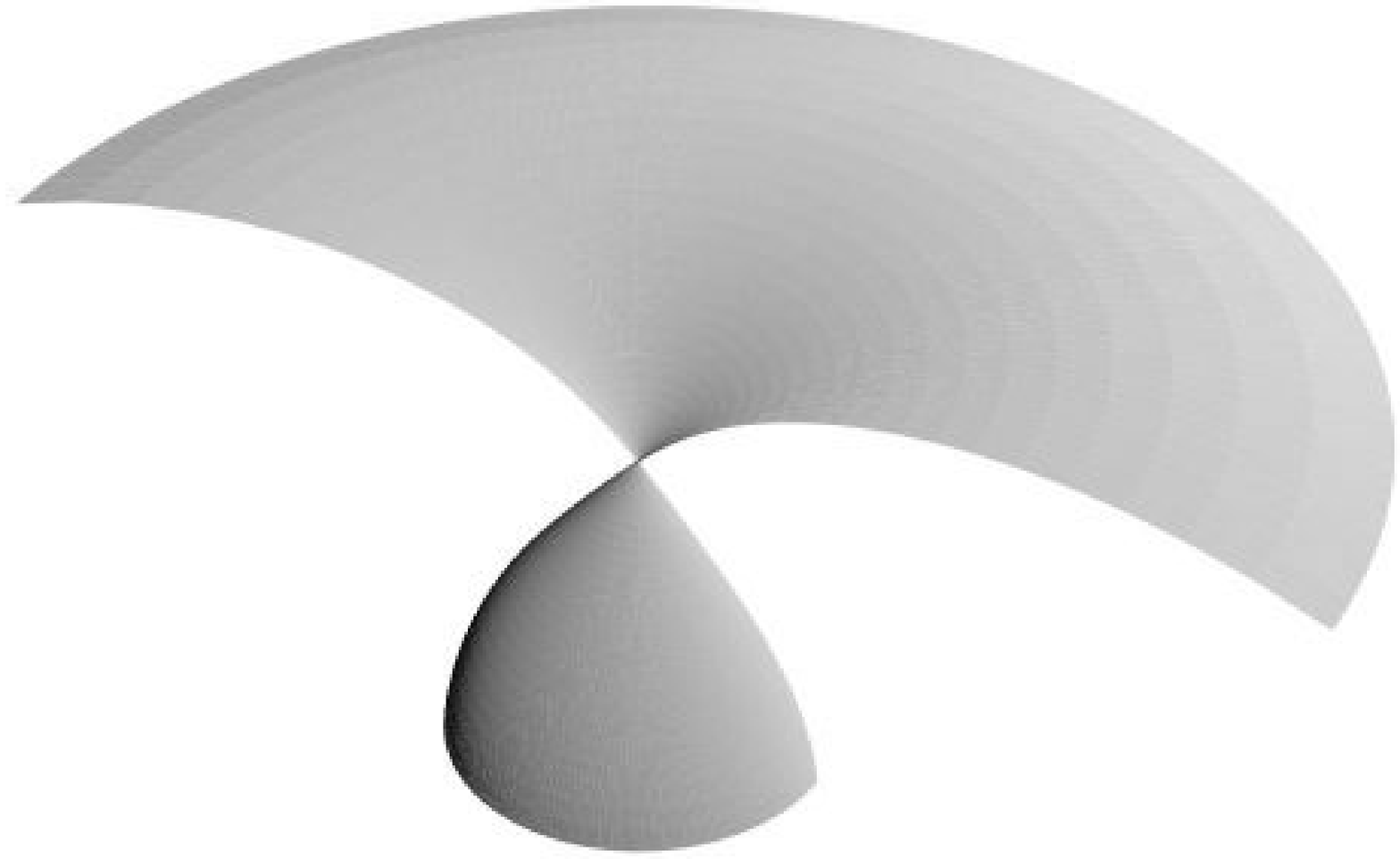} \\
 $\{z\in\mathbb{C}\,;\,
    \begin{subarray}{c}e^{-5}<|z|<e^5\\0<\arg z<\pi\end{subarray}\}.$ &
 $\{z\in\mathbb{C}\,;\,
    \begin{subarray}{c}e^{-5}<|z|<e^5\\\pi <\arg z<(3/2)\pi\end{subarray}\}.$ &
 $\{z\in\mathbb{C}\,;\,
    \begin{subarray}{c}e^{-2}<|z|<e^2\\0<\arg z<\pi\end{subarray}\}.$
\end{tabular}
\end{center}
\caption{Pictures of Example \ref{ex:cat}, where $\mu =1.2$.}
\end{figure} 

To produce further examples, we consider a relationship between CMC $1$ faces 
and CMC $1$ immersions in the hyperbolic space $\mathbb{H}^3$, and shall give 
a method for transferring from CMC $1$ immersions in $\mathbb{H}^3$ to 
CMC $1$ faces in $\mathbb{S}^3_1$. 

Let $\hat f:M=\overline{M}\setminus\{p_1,\dots ,p_n\}\to\mathbb{H}^3$ be a 
reducible CMC $1$ immersion, whose first fundamental form $d\hat s^2$ has 
finite total curvature and is complete, where we define $\hat f$ to be 
{\em reducible} if there exists a holomorphic null lift $F$ of 
$\hat f$ such that the image of the monodromy representation is in 
\[
U(1)=\left\{\left.\begin{pmatrix}e^{i\theta}&0\\0&e^{-i\theta}\end{pmatrix}
            \right|\theta\in\mathbb{R}\right\} .
\]
Let $(g,\omega)$ be the Weierstrass data associated to $F$, that is, 
$(g,\omega)$ satisfies 
\[
F^{-1}dF=\begin{pmatrix}g&-g^2\\1&-g\end{pmatrix}\omega .
\]
Then $|g|$ and $|\omega|$ are single-valued on $M$, 
as seen in the proof of Proposition \ref{pr:dshat2}. 
Assume that the absolute value of the secondary Gauss map is 
not equal to $1$ at all ends $p_1,\dots ,p_n$. 
Define $f:=Fe_3F^*$. Then $f$ is defined on $M$ as well. 
Furthermore, we have the following proposition:

\begin{proposition}\label{pr:UYmax5.4}
The CMC $1$ face $f:M\to\mathbb{S}^3_1$ defined as above, using $\hat f$ and 
its lift $F$ with monodromy in $U(1)$, is complete and of 
finite type with only elliptic ends. 
Moreover, an end of $\hat f$ is embedded if and only if the corresponding 
end of $f$ is embedded. 
\end{proposition}

\begin{proof}
Fix an end $p_j$ and assume $|g(p_j)|<1$. 
Then we can take a neighborhood $U_j$ such that $|g|^2<1-\varepsilon$ holds on 
$U_j$, where $\varepsilon\in (0,1)$ is a constant. In this case, 
\[
ds^2=(1-|g|^2)^2\omega\bar\omega\ge\varepsilon^2\omega\bar\omega
\ge\frac{\varepsilon^2}{2}(1+|g|^2)^2\omega\bar\omega
=\frac{\varepsilon^2}{2}d\hat s^2
\]
holds on $U_j$.  
Since $d\hat s^2$ is complete at $p_j$, $ds^2$ is also complete. 
Moreover, the Gaussian curvatures $K$ and $K_{d\hat s^2}$ satisfy 
\[
Kds^2=\frac{4dgd\bar g}{(1-|g|^2)^2}
\le\left(\frac{2}{\varepsilon}-1\right)^2\frac{4dgd\bar g}{(1+|g|^2)^2}
=\left(\frac{2}{\varepsilon}-1\right)^2(-K_{d\hat s^2})d\hat s^2.
\]
Hence $ds^2$ is of finite total curvature at the end $p_j$.

On the other hand, if $|g(p_j)|>1$, we can choose the neighborhood $U_j$ such 
that $|g|^{-2}<1-\varepsilon$ holds on $U_j$. Then 
\[
ds^2=(1-|g|^{-2})^2|g|^4\omega\bar\omega\ge\varepsilon^2|g|^4\omega\bar\omega
\ge\varepsilon^2(1+|g|^2)^2\omega\bar\omega=\varepsilon^2d\hat s^2.
\]
Hence $ds^2$ is complete at $p_j$. 
Moreover, since
\[
Kds^2=\frac{4dgd\bar g}{(1-|g|^2)^2}
\le\left(\frac{2}{\varepsilon}+1\right)^2\frac{4dgd\bar g}{(1+|g|^2)^2}
=\left(\frac{2}{\varepsilon}+1\right)^2(-K_{d\hat s^2})d\hat s^2,
\]
$ds^2$ is of finite total curvature.
The proof of the final sentence of the proposition follows from  the proof of 
Theorem \ref{lm:sub2}, by showing $m_1=1$ for both $f$ and $\hat f$. 
\end{proof}

Moreover, \cite[Theorem 3.3]{UY1} shows that for each 
$\lambda\in\mathbb{R}\setminus\{0\}$, $(\lambda g,\lambda^{-1}\omega)$ induces 
a CMC $1$ immersion $\hat f_\lambda :M\to\mathbb{H}^3$, where $(g,\omega)$ is
stated in Proposition \ref{pr:UYmax5.4}. 
Thus we have the following theorem:

\begin{theorem}\label{th:UYmax5.4}
Let $\hat f:M\to\mathbb{H}^3$ be a reducible complete CMC $1$ immersion of 
finite total curvature with $n$ ends. Then there exists the holomorphic null 
lift $F$ so that the image of the monodromy representation is in $U(1)$. 
Let $(g,\omega)$ be the Weierstrass data associated to $F$. 
Then there exist  $m$ $(0\le m\le n)$ positive real numbers 
$\lambda_1,\dots ,\lambda_m\in\mathbb{R}^+$ such that 
$f_\lambda :M\to\mathbb{S}^3_1$, induced from the Weierstrass data 
$(\lambda g,\lambda^{-1}\omega)$, is a complete CMC $1$ face of finite type 
with only elliptic ends for any 
$\lambda\in\mathbb{R}\setminus\{0,\pm\lambda_1,\dots ,\pm\lambda_m\}$. 
\end{theorem}

\begin{proof}
Let $\hat f$, $F$ and $(g,\omega)$ be as above and set 
$M=\overline{M}\setminus\{p_1,\dots ,p_n\}$.  
Then by \cite[Theorem 3.3]{UY1}, there exists a $1$-parameter family of 
reducible complete CMC $1$ immersions $\hat f_\lambda:M\to\mathbb{H}^3$ of 
finite total curvature with $n$ ends. 
Define $\lambda_j\in\mathbb{R}^+\cup\{0,\infty\}$ ($j=1,\dots ,m$) as 
\[
\lambda_j=
\left\{\begin{array}{ll}
       0        & \mbox{if $|g(p_j)|=\infty$}, \\
       \infty   & \mbox{if $|g(p_j)|=0$}, \\
  |g(p_j)|^{-1} & \mbox{otherwise}. 
       \end{array}\right.
\]
Then for any $\lambda\in 
(\mathbb{R}\cup\{\infty\})\setminus\{0,\pm\lambda_1,\dots ,\pm\lambda_m\}$, 
$|\lambda g(p_j)|\ne1$ for all $p_j$, and hence 
$f_\lambda :M\to\mathbb{S}^3_1$ induced from the Weierstrass data 
$(\lambda g,\lambda^{-1}\omega)$ is a complete CMC $1$ face of finite type 
with only elliptic ends.
\end{proof}

Complete CMC $1$ immersions with low total curvature and low dual 
total curvature in $\mathbb{H}^3$ were classified in \cite{RUY4, RUY5}. 
Applying Theorem \ref{th:UYmax5.4} to the reducible examples in their 
classification, we have the following:

\begin{corollary}\label{co:exmples}
There exist the following twelve types of complete CMC $1$ faces 
$f:M\to\mathbb{S}^3_1$ of finite type with elliptic ends:
\[
\begin{array}{llll}
{\bf O}(0),     & {\bf O}(-5),    & {\bf O}(-2,-3), & {\bf O}(-1,-1,-2), \\
{\bf O}(-4),    & {\bf O}(-6),    & {\bf O}(-2,-4), & {\bf O}(-1,-2,-2), \\
{\bf O}(-2,-2), & {\bf O}(-1,-4), & {\bf O}(-3,-3), & {\bf O}(-2,-2,-2),
\end{array}
\]
where $f$ is of {\em type} ${\bf O}(d_1,\dots,d_n)$ when 
$M=(\mathbb{C}\cup\{\infty\})\setminus\{p_1,\dots,p_n\}$ 
and $Q$ has order $d_j$ at each end $p_j$.  
\end{corollary}

Furthermore, reducible complete CMC $1$ immersions of genus zero with 
an arbitrary number of regular ends and one irregular end and finite total 
curvature are constructed in \cite{MU}, 
using an analogue of the so-called UP-iteration.  
Applying Theorem \ref{th:UYmax5.4} to their results, we have the following: 

\begin{corollary}\label{co:UPiteration}
Set $M=\mathbb{C}\setminus\{p_1,\dots,p_n\}$ for arbitrary $n\in\mathbb{N}$.  
Then there exist choices for $p_1,\dots,p_n$ so that there exist complete 
CMC $1$ faces $f:M\to\mathbb{S}^3_1$ of finite type 
with $n$ regular elliptic ends and one irregular elliptic end. 
\end{corollary}

\end{document}